\documentclass[final]{siamltex}

\usepackage[USenglish]{babel} %francais, polish, spanish, ...
\usepackage[T1]{fontenc}
\usepackage[ansinew]{inputenc}
\usepackage{lmodern}

\usepackage{graphicx} %%For loading graphic files

\usepackage{amsmath}
\usepackage{amsfonts}
\usepackage{amssymb}
\usepackage{url}

\usepackage{tikz}
\usetikzlibrary{scopes}
\usetikzlibrary{chains}
\usetikzlibrary{positioning}
\usetikzlibrary{arrows}

\usepackage[caption=false]{subfig}

\usepackage[sharp]{easylist}

\usepackage{todonotes}

\usepackage{ifpdf}
\ifpdf
  \usepackage{epstopdf} 
\fi

\usepackage[letterpaper]{geometry}
\newcommand{\etal}{\textit{et al.\ }}
\newcommand{\CC}{C\nolinebreak\hspace{-.05em}\raisebox{.4ex}{\tiny\bf +}\nolinebreak\hspace{-.10em}\raisebox{.4ex}{\tiny\bf +}}

\newtheorem{Example}[theorem]{Example}

\newtheorem{Definition}[theorem]{Definition}
\newtheorem{Remark}[theorem]{Remark}

\providecommand{\OO}[1]{\ensuremath{\mathcal{O}\bigl(#1\bigr)}}
\usepackage[colorlinks=true, citecolor=blue, filecolor=black,linkcolor=red, urlcolor=blue, hypertexnames=false, pdftex]{hyperref}

\title{Adaptive Smolyak Pseudospectral Approximations}
\author{Patrick R.\ Conrad and Youssef M.\ Marzouk}

\begin{document}

\maketitle

\begin{abstract}
  Polynomial approximations of computationally intensive models are
  central to uncertainty quantification. This paper describes an
  adaptive method for non-intrusive pseudospectral approximation, based on Smolyak's
  algorithm with generalized sparse grids. We 
  rigorously analyze and extend 
  the non-adaptive method proposed in
  \cite{Constantine2012}, and compare it to a common alternative
  approach for using sparse grids to construct polynomial
  approximations, direct quadrature. Analysis of direct quadrature
  shows that $\mathcal{O}(1)$ errors are an intrinsic property of some configurations of the method, as a
  consequence of internal aliasing. We provide precise conditions,
  based on the chosen polynomial basis and quadrature rules, under
  which this aliasing error occurs. 
  We then establish theoretical results on the accuracy of Smolyak
  pseudospectral approximation, and show that the Smolyak
  approximation avoids internal aliasing and makes far more effective
  use of sparse function evaluations. These results are applicable to
  broad choices of quadrature rule and generalized sparse
  grids. Exploiting this flexibility, we introduce a greedy heuristic
  for adaptive refinement of the pseudospectral approximation. We
  numerically demonstrate convergence of the algorithm on the Genz
  test functions, and illustrate the accuracy and efficiency of the
  adaptive approach on a realistic chemical kinetics problem.

  % Polynomial approximations of computationally intensive models are
  % central to uncertainty quantification. This paper describes an
  % adaptive method for pseudospectral approximation, based on Smolyak's
  % algorithm with generalized sparse grids. We rigorously develop the
  % method and compare it to a common alternative approach for using
  % sparse grids to construct polynomial approximations, direct
  % quadrature. Analysis of direct quadrature shows that the large
  % errors observed in \cite{Constantine2012} are an intrinsic risk of
  % the method, a consequence of internal aliasing; we provide precise
  % conditions, based on the chosen polynomial basis and quadrature
  % rules, under which this aliasing error occurs. 
  % %
  % We then establish [key!] theoretical results on the accuracy of Smolyak
  % pseudospectral approximation, and show that the Smolyak
  % approximation avoids internal aliasing and makes far more effective
  % use of sparse function evaluations. These results are applicable to
  % broad choices of quadrature rule and generalized sparse grids.
  % Exploiting this flexibility, we introduce a greedy heuristic for
  % adaptive refinement of the pseudospectral approximation. We
  % numerically demonstrate convergence of the algorithm on the Genz
  % test functions, and illustrate the accuracy and efficiency of the
  % adaptive approach on a realistic chemical kinetics problem.
  
\end{abstract}

\begin{keywords} 
  Smolyak algorithms, sparse grids, orthogonal polynomials,
  pseudospectral approximation, approximation
  theory, uncertainty quantification
\end{keywords}

\begin{AMS}
41A10, 41A63, 47A80, 65D15, 65D32
\end{AMS}

\pagestyle{myheadings}
\thispagestyle{plain}
\markboth{CONRAD AND MARZOUK}{\MakeUppercase{Adaptive Smolyak Pseudospectral Approximations}}

%% Chapters %%%%%%%%%%%%%%%%%%%%%%%%%%%%%%%%%%%%%%%%%%%%%%%%%
%% ==> Write your text here or include other files.

\section{Introduction}

A central issue in the field of uncertainty quantification is understanding the
response of a model to random inputs. When model evaluations are
computationally intensive, techniques for \textit{approximating} the
model response in an efficient manner are essential. Approximations
may be used to evaluate moments or the probability distribution of a
model's outputs, or to evaluate sensitivities of model outputs with
respect to the inputs~\cite{LeMaitre2010, Xiu2010,
  Sudret2008}. Approximations may also be viewed as \textit{surrogate
  models} to be used in optimization~\cite{march:cmo:2012} or
inference~\cite{Marzouk2007}, replacing the full model entirely.

Often one is faced with black box models that can only be evaluated at
designated input points. We will focus on constructing multivariate
polynomial approximations of the input-output relationship generated
by such a model; these approximations offer fast convergence for
smooth functions and are widely used. One common strategy for
constructing a polynomial approximation is interpolation, where
interpolants are conveniently represented in Lagrange form
\cite{Babuska2007,Xiu2005}.  Another strategy is
projection, particularly orthogonal projection with respect to some
inner product. The results of such a projection are conveniently
represented with the corresponding family of orthogonal polynomials
\cite{Canuto2006,LeMaitre2010,Xiu2002}. When the inner product is
chosen according to the input probability measure, this construction
is known as the (finite dimensional) polynomial chaos expansion (PCE)
\cite{Ghanem1991,Soize2004,ernst:ocg:2012}. Interpolation
and projection are closely linked, particularly when projection is
computed via discrete model evaluations. Moreover, one can always
realize a change of basis \cite{gander:cbp:2005} for the polynomial
resulting from either operation. Here we will favor orthogonal
polynomial representations, as they are easy to manipulate and their
coefficients have a useful interpretation in probabilistic settings.

This paper discusses \emph{adaptive Smolyak pseudospectral
  approximation}, an accurate and computationally efficient 
approach to constructing multivariate polynomial chaos expansions.
Pseudospectral methods allow the construction of polynomial
approximations from point evaluations of a function \cite{Canuto2006,
  Boyd2001}. We combine these methods with \emph{Smolyak's algorithm},
a general strategy for sparse approximation of linear operators on
tensor product spaces, which saves computational effort by weakening
the assumed coupling between the input dimensions.  Gerstner \&
Griebel~\cite{Gerstner2003} and Hegland~\cite{Hegland2003} developed adaptive variants of Smolyak's
algorithm for numerical integration and illustrated the effectiveness
of on-the-fly heuristic adaptation. We extend their approach to
the pseudospectral approximation of functions.  Adaptivity is expected
to yield substantial efficiency gains in high
dimensions---particularly for functions with anisotropic dependence on
input parameters and functions whose inputs might not be strongly
coupled at high order.

%%%%%%%%%%%%%%%%%%%%%%%%%%

Previous attempts to extend pseudospectral methods to multivariate
polynomial approximation with sparse model evaluations employed
{ad hoc} approaches that are not always accurate. A common
procedure has been to use sparse quadrature, or even
dimension-adaptive sparse quadrature, to evaluate polynomial
coefficients directly \cite{Xiu2007,LeMaitre2010}.
This leads to at least two difficulties. First, the truncation of the
polynomial expansion must be specified independently of the quadrature
grid, yet it is unclear how to do this, particularly for anisotropic
and generalized sparse grids. Second, unless one uses excessively
high-order quadrature, significant aliasing errors may
result. Constantine \etal \cite{Constantine2012} provided the first
clear demonstration of these aliasing errors and proposed a Smolyak
algorithm that does not share them. That work also demonstrated a link
between Smolyak pseudospectral approximation and an extension to
Lagrange interpolation called \emph{sparse interpolation}, which uses
function evaluations on a sparse grid and has well characterized
convergence properties \cite{Nobile2007, Barthelmann2000}.

The first half of this work performs a theoretical analysis, placing
the solution from \cite{Constantine2012} in the broader context of
Smolyak constructions, and explaining the origin of the observed
aliasing errors for general (e.g., anisotropic) choices of sparse grid
and quadrature rule. We do so by using the notion of polynomial
exactness, without appealing to interpolation properties of particular
quadrature rules. We establish conditions under which tensorized
approximation operators are exact for particular polynomial inputs,
then apply this analysis to the specific cases of quadrature and
pseudospectral approximation; these cases are closely related and
facilitate comparisons between Smolyak pseudospectral
algorithms and direct quadrature.
%
% We accomplish this by studying when computable approximations are
% exact for polynomial inputs, establishing useful properties at finite
% orders. We apply this analysis to the general case of tensor
% approximation operators, in addition to the specific cases of
% quadrature and pseudospectral approximation; these cases are closely
% related and facilitate the desired comparisons between Smolyak
% pseudospectral algorithms and direct quadrature.
%
Section \ref{sec:approximations} develops \textit{computable}
one-dimensional and tensorized approximations for these
settings. Section \ref{sec:smolyak} describes general Smolyak
algorithms and their properties, yielding our principal theorem about
the polynomial exactness of Smolyak approximations, and then applies
these results to quadrature and pseudospectral approximation. Section
\ref{sec:comparison} compares the Smolyak approach to conventional
direct quadrature. Our error analysis of direct quadrature shows why
the approach goes wrong and allows us to draw an important conclusion:
in almost all cases, direct quadrature is not an appropriate method
for constructing polynomial expansions and should be superseded by
Smolyak pseudospectral methods.

These results provide a rigorous foundation for \textit{adaptivity},
which is the second focus of this paper. Adaptivity makes it possible
to harness the full flexibility of Smolyak algorithms in practical
settings. Section \ref{sec:adaptive} introduces a fully adaptive
algorithm for Smolyak pseudospectral approximation, which uses a
single tolerance parameter to drive iterative refinement of both the
polynomial approximation space and the corresponding collection of
model evaluation points. As the adaptive method is largely heuristic,
Section \ref{sec:experiments} demonstrates the benefits of this
approach with numerical examples.

\section{Full tensor approximations}
\label{sec:approximations}

Tensorization is a common approach for lifting one-dimensional operators to higher dimensions. Not only are tensor products computationally convenient, but they provide much useful structure for analysis. In this section, we develop some essential background for computable tensor approximations, then apply it to problems of (i) approximating integrals with numerical quadrature; and (ii) approximating projection onto polynomial spaces with pseudospectral methods. In particular, we are interested in analyzing the errors associated with these approximations and in establishing conditions under which the approximations are \textit{exact}.

\subsection{General setting}
Consider a collection of one-dimensional linear operators $\mathcal{L}^{(i)}$, where $(i)$ indexes the operators used in different dimensions. In this work, $\mathcal{L}^{(i)}$ will be either an integral operator or an orthogonal projector onto some polynomial space. We can extend a collection of these operators into higher dimensions by constructing the tensor product operator
\begin{equation}
\label{e:tpoperator}
\mathcal{L}^{(\mathbf{d})} := \mathcal{L}^{(1)} \otimes \cdots \otimes \mathcal{L}^{(d)} .
\end{equation}
The one-dimensional operators need not be identical; the properties of the resulting tensor operator are constructed independently from each dimension. The bold parenthetical superscript refers to the tensor operator instead of the constituent one-dimensional operators.

As the true operators are not available computationally, we work with a convergent sequence of computable approximations, $\mathcal{L}^{(i)}_m$, such that 
\begin{equation}
\| \mathcal{L}^{(i)} - \mathcal{L}^{(i)}_m \| \to 0 \ \mathrm{ as } \ m \to \infty
\end{equation}
in some appropriate norm. Taking the tensor product of these approximations provides an approximation to the full tensor operator, $\mathcal{L}^{(\mathbf{d})}_\mathbf{m}$, where the level of the approximation may be individually selected in each dimension, so the tensor approximation is identified by a multi-index $\mathbf{m}$. Typically, and in the cases of interest in this work, the tensor approximation will converge in the same sense as the one-dimensional approximation as all components of $\mathbf{m} \to \infty$. 

An important property of approximation algorithms is whether they are \textit{exact} for some inputs; characterizing this set of inputs allows us to make useful statements at finite order.

\begin{Definition}[Exact Sets]
For an operator $\mathcal{L}$ and a corresponding approximation $\mathcal{L}_m$, define the exact set as $\mathcal{E}(\mathcal{L}_m) := \{f: \mathcal{L}(f) = \mathcal{L}_m(f)\}$ and the half exact set $\mathcal{E}_2(\mathcal{L}_m) := \{f : \mathcal{L}(f^2) = \mathcal{L}_m(f^2)\}$.
\end{Definition}

The half exact set will help connect the exactness of a quadrature rule to that of the closely related pseudospectral operators. This notation is useful in proving the following lemma, which relates the exact sets of one-dimensional approximations and tensor approximations.

\begin{lemma}
\label{thm:tensorAccuracy}
If a tensor approximation $\mathcal{L}^{(\mathbf{d})}_\mathbf{m}$ is constructed from one-dimensional approximations $\mathcal{L}^{(i)}_m$ with known exact sets, then

\begin{equation}
\mathcal{E}(\mathcal{L}^{(1)}_{{m}_1}) \otimes \cdots \otimes \mathcal{E}(\mathcal{L}^{(d)}_{{m}_d}) \,  \subseteq \, \mathcal{E}(\mathcal{L}^{(\mathbf{d})}_\mathbf{m})
\end{equation}
\end{lemma}

\proof{It is sufficient to show that the approximation is exact for an arbitrary monomial input $f (\mathbf{x} ) = f^{(1)}(x^{(1)}) f^{(2)}(x^{(2)}) \cdots f^{(d)}(x^{(d)})$ where $f^{(i)}(x^{(i)}) \in \mathcal{E}(\mathcal{L}^{(i)}_{{m}_i})$, because we may extend to sums by linearity:
\begin{eqnarray*}
\mathcal{L}^{(\mathbf{d})}_{\mathbf{m}}(f^{(1)} \cdots f^{(d)}) & = & \mathcal{L}^{(1)}_{{m}_1}(f^{(1)}) \otimes \cdots \otimes \mathcal{L}^{(d)}_{{m}_d}(f^{(d)})  \\
& = & \mathcal{L}^{(1)}(f^{(1)}) \otimes \cdots \otimes \mathcal{L}^{(d)}(f^{(d)}) = \mathcal{L}^{(\mathbf{d})} ( f ) .
\end{eqnarray*}
The first step uses the tensor product structure of the operator and the second uses the definition of exact sets.
}\endproof

\subsection{Multi-indices}
Before continuing, we must make a short diversion to multi-indices, which provide helpful notation when dealing with tensor problems. A multi-index is a vector $\mathbf{i} \in \mathbb{N}^d_0$. An important notion for multi-indices is that of a  \textit{neighborhood}.

\begin{Definition}[Neighborhoods of multi-indices]
A \emph{forward neighborhood} of a multi-index $\mathbf{k}$ is the multi-index set $n_f(\mathbf{k}) := \{\mathbf{k}+\mathbf{e}_i: \forall i \in \{1 \ldots d\} \}$, where $\mathbf{e}_i$ are the canonical unit vectors. The \emph{backward neighborhood} of a multi-index $\mathbf{k}$ is the multi-index set $n_b(\mathbf{k}) := \{\mathbf{k}-\mathbf{e}_i : \forall i \in \{1 \ldots d\}, \mathbf{k}-\mathbf{e}_i \in \mathbb{N}^d_0  \}$. 
\end{Definition}

Smolyak algorithms rely on multi-index sets that are \emph{admissible}.

\begin{Definition}[Admissible multi-indices and multi-index sets]
A multi-index $\mathbf{k}$ is admissible with respect to a multi-index set $\mathcal{K}$ if $n_b(\mathbf{k}) \subseteq \mathcal{K}$. A multi-index set $\mathcal{K}$ is admissible if every $\mathbf{k} \in \mathcal{K}$ is admissible with respect to $\mathcal{K}$.
\end{Definition}

Two common admissible multi-index sets with simple geometric structure are \emph{total order} multi-index sets and \emph{full tensor} multi-index sets. One often encounters total order sets in the sparse grids literature and full tensor sets when dealing with tensor grids of points. The total order multi-index set $\mathcal{K}^{t}_n$ comprises those multi-indices that lie within a $d$-dimensional simplex of side length $n$:
\begin{equation}
\mathcal{K}^{t}_n := \{\mathbf{k} \in \mathbb{N}_0^d: \|\mathbf{k}\|_1 \leq n\}
\end{equation}
The full tensor multi-index set $\mathcal{K}^{f}_\mathbf{n}$ is the complete grid of indices bounded term-wise by a multi-index $\mathbf{n}$:
\begin{equation}
\mathcal{K}^{f}_\mathbf{n} := \{\mathbf{k} \in \mathbb{N}_0^d: \forall i \in \{1 \ldots d\}, \, {k}_i < {n}_i\}
\end{equation}

\subsection{Integrals and quadrature}

Let $X^{(i)}$ be an open or closed interval of the real line $\mathbb{R}$. Then we define the weighted integral operator in one dimension as follows:
\begin{equation}
\mathcal{I}^{(i)}(f) := \int_{X^{(i)}} f(x) w^{(i)}(x)\, dx
\end{equation}
where $f:X^{(i)} \to \mathbb{R}$ is some real-valued function and $w^{(i)}: X^{(i)} \to \mathbb{R}^+$ is an integrable weight function. We may extend to higher dimensions by forming the tensor product integral $\mathcal{I}^{(\mathbf{d})}$, which uses separable weight functions and Cartesian product domains. 

Numerical quadrature approximates the action of an integral operator $\mathcal{I}^{(i)}$ with a weighted sum of point evaluations. For some family of quadrature rules, we write the ``level $m$'' quadrature rule, comprised of $p^{(i)}(m): \mathbb{N} \rightarrow \mathbb{N}$ points, as 
\begin{equation}
\mathcal{I}^{(i)}(f) \approx \mathcal{Q}^{(i)}_m(f) := \sum_{j=1}^{p^{(i)}(m)} w_j^{(i)} f(x_j^{(i)}) .
\end{equation}
We call $p^{(i)}(m)$ the growth rate of the quadrature rule, and its form depends on the quadrature family; some rules only exist for certain numbers of points and others may be tailored, for example, to produce linear or exponential growth in the number of quadrature points with respect to the level.  

\begin{figure}
	\centering
		\includegraphics[scale=.7]{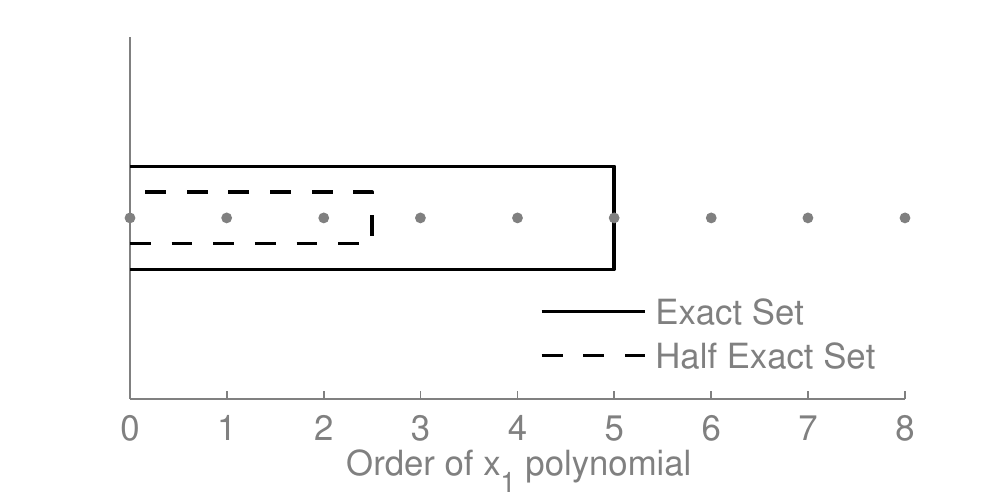}
	\caption{Consider, one-dimensional Gaussian quadrature rule with three points, $\mathcal{Q}^{(i)}_3$, which is exact for fifth degree polynomials. This diagram depicts the exact set, $\mathcal{E}(\mathcal{Q}^1_3)$, and half exact set, $\mathcal{E}_2(\mathcal{Q}^1_3)$, of this quadrature rule.}
	\label{fig:Quadrature1D}
\end{figure}

Many quadrature families are exact if $f$ is a polynomial of a degree $a^{(i)}(m)$ or less, which allows us to specify a well-structured portion of the exact set for these quadrature rules:

\begin{align}
\mathbb{P}_{a^{(i)}(m)} & \subseteq \mathcal{E}(\mathcal{Q}^{(i)}_m)\\
\mathbb{P}_{\mathrm{floor}(a^{(i)}(m)/2)} & \subseteq \mathcal{E}_2(\mathcal{Q}^{(i)}_m),
\end{align}
where $\mathbb{P}_a$ is the space of polynomials of degree $a$ or less. It is intuitive and useful to draw the exact set as in Figure \ref{fig:Quadrature1D}. For this work, we rely on quadrature rules that exhibit polynomial accuracy of increasing order, which is sufficient to demonstrate convergence for functions in $L^2$ \cite{Canuto2006}.

\begin{figure}
	\centering
		\includegraphics[scale=.6]{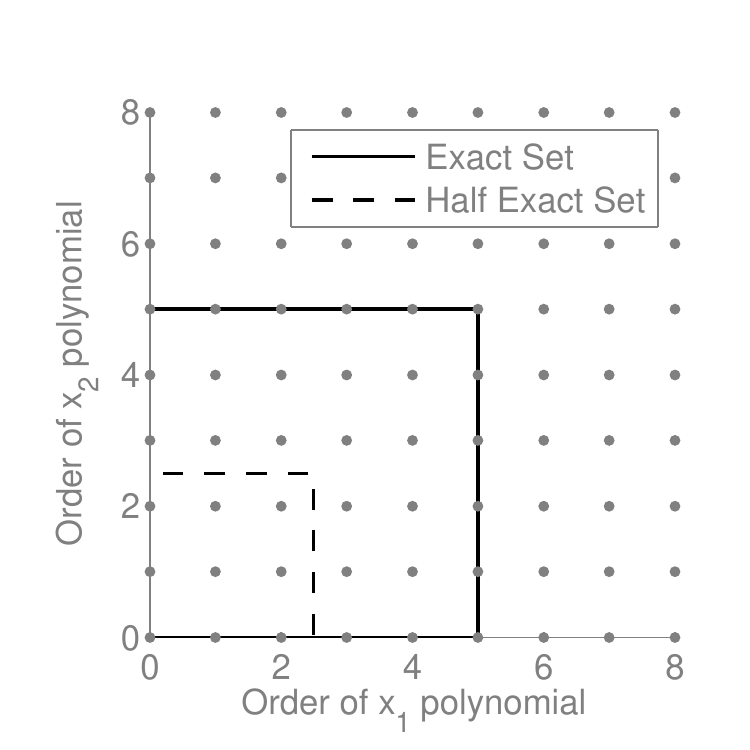}
	\caption{Consider a two dimensional quadrature rule constructed from three point Gaussian quadrature rules, $\mathcal{Q}^{(\mathbf{2})}_{(3,3)}$. This diagram depicts the exact set, $\mathcal{E}(\mathcal{Q}^{(\mathbf{2})}_{(3,3)})$, and half exact set, $\mathcal{E}_2(\mathcal{Q}^{(\mathbf{2})}_{(3,3)})$.}
	\label{fig:SingleTensorQuadrature}
\end{figure}

Tensor product quadrature rules are straightforward approximations of tensor product integrals that inherit convergence properties from the one-dimensional case. The exact set of a tensor product quadrature rule includes the tensor product of the constituent approximations' exact sets, as guaranteed by Lemma \ref{thm:tensorAccuracy} and depicted in Figure \ref{fig:SingleTensorQuadrature}.

%We can identify another portion of the exact set that arises from the tensor structure of the operator, which we use later in our analysis of sparse problems.
%
%\begin{Remark}
%Let $f$ be separable, without loss of generality, with respect to its first input, i.e., $f(\mathbf{x}) = f^{\prime}(x_1)f^{\prime \prime}(x_2, \ldots x_d)$. Further, let $\mathcal{I}^1(f^{\prime})=0$ and let $f^{\prime}(x_1) \in \mathcal{E}(\mathcal{Q}^{(1)}_{m_1})$. Then $\mathcal{Q}^{(\mathbf{d})}_\mathbf{m}(f) = \mathcal{I}^{(\mathbf{d})}(f) = 0$ for any $f^{\prime \prime}$. Hence, $f \in \mathcal{E}(\mathcal{Q}^{(\mathbf{d})}_\mathbf{m})$.
%\end{Remark}

\subsection{Polynomial projection}

A polynomial chaos expansion approximates a function with a weighted sum of orthonormal polynomials \cite{Wiener1938,Xiu2002}. Let $\mathcal{H}^{(i)} := L^2 \left (X^{(i)}, w^{(i)} \right )$ be the separable Hilbert space of square-integrable functions $f: X^{(i)} \rightarrow \mathbb{R}$, with inner product defined by the weighted integral $\langle f,g \rangle = \mathcal{I}^{(i)}(fg)$, and $w^{(i)}(x)$ normalized so that it may represent a probability density. Let $\mathbb{P}^{(i)}_n$ be the space of univariate polynomials of degree $n$ or less. Now let $\mathcal{P}_n^{(i)}:\mathcal{H}^{(i)} \rightarrow \mathbb{P}^{(i)}_n$ be an orthogonal projector onto this subspace, written in terms of polynomials $\{ \psi_j^{(i)}(x): j \in \mathbb{N}_0 \}$ orthonormal with respect to the inner product of $\mathcal{H}^{i}$:
\begin{equation}
\label{eq:truncated1DPCE}
\mathcal{P}^{(i)}_n(f) := \sum_{j=0}^n \left \langle f(x), \psi^{(i)}_j(x) \right \rangle \psi^{(i)}_j(x) = \sum_{i=0}^n f_j \psi^{(i)}_j(x).
\end{equation}
The polynomial space $\mathbb{P}^{(i)}_n$ is of course the \textit{range} of the projection operator.
These polynomials are dense in $\mathcal{H}^{(i)}$, so the polynomial approximation of any $f \in \mathcal{H}^{(i)}$ converges in the $L^2$ sense as $n \to \infty$ \cite{Canuto2006, Xiu2002}. 
If $f \in \mathcal{H}^{(i)}$, the coefficients must satisfy $\sum_{i=0}^\infty f_i^2 < \infty$. 

Projections with finite degree $n$ omit terms of the infinite series, thus incurring \emph{truncation error}. We can write this error as 
\begin{equation}
\label{eq:truncation}
\left \| f - \mathcal{P}^{(i)}_n(f) \right \|_2^2 = \left \| \sum_{j=n+1}^\infty f_j \psi_j^{(i)} \right \|_2^2 = \sum_{j=n+1}^\infty f_j^2 < \infty.
\end{equation}
Hence, we may reduce the truncation error to any desired level by increasing $n$, removing terms from the sum in (\ref{eq:truncation}) \cite{Canuto2006,Hesthaven2007}.

The $d$-dimensional version of this problem requires approximating functions in the Hilbert space $\mathcal{H}^{(\mathbf{d})} := \mathcal{H}^{(1)} \otimes \cdots \otimes \mathcal{H}^{(d)}$ via a tensor product basis of the univariate polynomials defined above:
\begin{equation}
\mathcal{P}_\mathbf{n}^{(\mathbf{d})}(f) = \sum_{{i}_1=0}^{{n}_1} \ldots \sum_{{i}_d=0}^{{n}_d} \left \langle f \Psi_\mathbf{i} \right \rangle \Psi_\mathbf{i}
\end{equation}
where $\Psi_\mathbf{i}(\mathbf{x}) := \prod_{j=1}^d \psi_{i_j}^{(j)} \left ( x^{(j)} \right )$. The multi-index $\mathbf{n}$ tailors the range of the projection to include a rectangular subset of polynomials.

As in the one-dimensional case, truncation induces error equal to the sum of the squares of the omitted coefficients, which we may similarly reduce to zero as ${n}_i \to \infty$, $\forall i$. The multivariate polynomial expansion also converges in an $L^2$ sense for any $f \in \mathcal{H}^{(\mathbf{d})}$ \cite{Canuto2006}

\subsection{Aliasing errors in pseudospectral approximation}
\label{sec:aliasing}
The inner products defining the expansion coefficients above are not directly computable. Pseudospectral approximation provides a practical non-intrusive algorithm by approximating these inner products with quadrature. Define the pseudospectral approximation in one dimension as
\begin{eqnarray}
\mathcal{S}_m^{(i)} (f) &:= &\sum_{j=0}^{q^{(i)}(m)} \mathcal{Q}^{(i)}_m \left ( f\psi_j^{(i)} \right ) \psi_j^{(i)}(x) \nonumber \\
&= & \sum_{j=0}^{q^{(i)}(m)} \tilde{f}_j \psi^{(i)}_j(x)
\end{eqnarray}
where $q^{(i)}(m)$ is the polynomial truncation at level $m$, to be specified shortly \cite{Canuto2006,Hesthaven2007}. Pseudospectral approximations are constructed around a level $m$ quadrature rule, and are designed to include as many terms in the sum as possible while maintaining accuracy. Assuming that $f \in L^2$, we can compute the $L^2$ error between the pseudospectral approximation and an exact projection onto the same polynomial space:
\begin{align}
\left \| \mathcal{P}_{q^{(i)}(m)}^{(i)}(f) - \mathcal{S}^{(i)}_m(f) \right \|_2^2  &= \left \| \sum_{j=0}^{q^{(i)}(m)} f_j \psi_j^{(i)} - \sum_{k=0}^{q^{(i)}(m)} \tilde{f}_k \psi_k^{(1)} \right \|^2_2 = \sum_{j=0}^{q^{(i)}(m)} (f_j - \tilde{f}_j )^2 
\end{align}
This quantity is the \emph{aliasing error} \cite{Canuto2006,Hesthaven2007}. The error is non-zero because quadrature in general only approximates integrals; hence each $\tilde{f}_i$ is an approximation of $f_i$. The pseudospectral operator also incurs truncation error, as before, which is orthogonal to the aliasing error. We can expand each approximate coefficient as
\begin{eqnarray}
\tilde{f}_j & = & \mathcal{Q}^{(i)}_m \left (f \psi^{(i)}_j \right ) \nonumber \\ 
	&= & \sum_{k=0}^\infty f_k \mathcal{Q}^{(i)}_m \left (\psi^{(i)}_j \psi^{(i)}_k \right) \nonumber \\
	&= &\sum_{k=0}^{q^{(i)}(m)} f_k \mathcal{Q}^{(i)}_m \left (\psi^{(i)}_j \psi^{(i)}_k \right ) + \sum_{l=q^{(i)}(m)+1}^\infty f_l \mathcal{Q}^{(i)}_m \left (\psi^{(i)}_j \psi^{(i)}_l \right ) .
\label{e:coeff}
\end{eqnarray}
The first step substitutes in the polynomial expansion of $f$, which we assume is convergent, and rearranges using linearity. The second step partitions the sum around the truncation of the pseudospectral expansion. Although the basis functions are orthonormal, $\langle  \psi^{(i)}_j, \psi^{(i)}_k \rangle = \delta_{jk}$, we cannot assume in general that the approximation $\mathcal{Q}^{(i)}_m \left  (\psi^{(i)}_j \psi^{(i)}_k \right ) = \delta_{jk}$. Now substitute (\ref{e:coeff}) back into the aliasing error expression:
\begin{equation}
\label{eqn:errorDecomp1D}
\sum_{j=0}^{q^{(i)}(m)} (f_j - \tilde{f}_j )^2  =  \sum_{j=0}^{q^{(i)}(m)} \left( f_j - \sum_{k=0}^{q^{(i)}(m)} f_k \, \mathcal{Q}^{(i)}_m \left (\psi^{(i)}_j \psi^{(i)}_k \right ) - \sum_{l=q^{(i)}(m)+1}^\infty f_l \, \mathcal{Q}^{(i)}_m \left (\psi^{(i)}_j \psi^{(i)}_l \right ) \right) ^2 
\end{equation}

This form reveals the intimate link between the accuracy of pseudospectral approximations and the polynomial accuracy of quadrature rules. All aliasing is attributed to the inability of the quadrature rule to determine the orthonormality of the basis polynomials, causing one coefficient to corrupt another. The contribution of the first two parenthetical terms on the right of (\ref{eqn:errorDecomp1D}) is called \emph{internal aliasing}, while the third term is called \emph{external aliasing}. Internal aliasing is due to inaccuracies in $\mathcal{Q}(\psi^{(i)}_j \psi^{(i)}_k)$ when \textit{both} $\psi^{(i)}_j$ and $\psi^{(i)}_k$ are included in the expansion, while external aliasing occurs when \textit{only one} of these polynomials is included in the expansion. For many practical quadrature rules (and for all those used in this work), if $j\neq k$ and $\psi^{(i)}_j \psi^{(i)}_k \notin \mathcal{E}(\mathcal{Q})$, and hence the discrete inner product is not zero, then $\| \mathcal{Q}(\psi^{(i)}_j \psi^{(i)}_k) \|_2 $ is \OO{1} \cite{Trefethen2008}. As a result, the magnitude of an aliasing error typically corresponds to the magnitude of the aliased coefficients.

In principle, both types of aliasing error are driven to zero by sufficiently powerful quadrature, but we are left to select $q^{(i)}(m)$ for a particular quadrature level $m$. External aliasing is \OO{1} in the magnitude of the truncation error, and thus it is driven to zero as long as $q^{(i)}(m)$ increases with $m$. Internal aliasing could be \OO{1} with respect to the function of interest, meaning that the procedure neither converges nor provides a useful approximation. Therefore, the obvious option is to include as many terms as possible while setting the internal aliasing to zero.

For quadrature rules with polynomial exactness, we may accomplish this by setting $q^{(i)}(m) = \mathrm{floor}( a^{(i)}(m)/2)$. This ensures that the internal aliasing of $\mathcal{S}^{(i)}_m$ is zero, because $\forall j,k \leq q^{(i)}(m),\\ \psi^{(i)}_j \psi^{(i)}_k \in \mathcal{E}(\mathcal{Q}^{(i)}_m)$. Equivalently, a pseudospectral operator $\mathcal{S}^{(i)}_m$ using quadrature $\mathcal{Q}^{(i)}_m$ has a range corresponding to the half exact set $\mathcal{E}_2(\mathcal{Q}^{(i)}_m)$. Alternatively, we may justify this choice by noting that it makes the pseudospectral approximation exact on its range, $\mathbb{P}_{q^{(i)}(m)}^{(i)} \subseteq \mathcal{E}(\mathcal{S}^{(i)}_m)$.

Given this choice of $q^{(i)}(m)$, we wish to show that the pseudospectral approximation converges to the true function, where the magnitude of the error is as follows:
\begin{equation}
\left \| f - \mathcal{S}^{(i)}_m(f) \right \|^2_2 = \sum_{j=0}^{q^{(i)}(m)}  \left( \sum_{k=q^{(i)}(m)+1}^\infty f_k \mathcal{Q} \left (\psi^{(i)}_j \psi^{(i)}_k \right ) \right)^2 + \sum_{l={q^{(i)}(m)} +1}^\infty f_l^2 .
\end{equation}
The two terms on right hand side comprise the external aliasing and the truncation error, respectively. We already know that the truncation error goes to zero as $q^{(i)}(m) \to \infty$. The external aliasing also vanishes for functions $f \in L^2$, as the truncated portion of $f$ likewise decreases \cite{Trefethen2008}. In the case of Gaussian quadrature rules, a link to interpolation provides precise rates for the convergence of the pseudospectral operator based on the regularity of $f$ \cite{Canuto2006}.

As with quadrature algorithms, our analysis of pseudospectral approximation in one dimension is directly extensible to multiple dimensions via full tensor products. We may thus conclude that $\mathcal{S}^{(\mathbf{d})}_\mathbf{m}$ converges to the projection onto the tensor product polynomial space in the same sense. The exact set follows Lemma $\ref{thm:tensorAccuracy}$, and hence the tensor product approximation inherits zero internal aliasing if suitable one-dimensional operators are used. 

\section{Smolyak algorithms}
\label{sec:smolyak}

Thus far, we have developed polynomial approximations of multivariate functions by taking tensor products of one-dimensional pseudospectral operators. Smolyak algorithms avoid the exponential cost of full tensor products when the input dimensions are not fully coupled, allowing the use of a telescoping sum to blend different lower-order full tensor approximations.

\begin{Example}
Suppose that $f(x,y) = x^{7} + y^7+x^3 y$. To construct a polynomial expansion with both the $x^7$ and $y^7$ terms, a full tensor pseudospectral algorithm would estimate all the polynomial terms up to $x^7y^7$, because tensor algorithms fully couple the dimensions. This mixed term is costly, requiring, in this case, an $8\times 8$ point grid for Gaussian quadratures. The individual terms can be had much more cheaply, using $8\times 1$, $1\times 8$, and $4 \times 2$ grids, respectively. Smolyak algorithms help realize such savings in practice.
\end{Example}

This section reviews the construction of Smolyak algorithms and presents a new theorem about the exactness of Smolyak algorithms built around arbitrary admissible index sets. We apply these results to quadrature and pseudospectral approximation, allowing a precise characterization of their errors.

\subsection{General Smolyak algorithms}
\label{s:generalsmolyak}
As in Section~\ref{sec:approximations}, assume that we have for every dimension $i=1\ldots d$ a convergent sequence $\mathcal{L}^{(i)}_k$ of approximations. Let $\mathcal{L}$ denote the collection of these sequences over all the dimensions. Define the difference operators
\begin{eqnarray}
\Delta^{(i)}_0 & := & \mathcal{L}^{(i)}_0 = 0, \\
\Delta^{(i)}_n & := & \mathcal{L}^{(i)}_n - \mathcal{L}^{(i)}_{n-1} .
\end{eqnarray}
For any $i$, we may write the exact or ``true'' operator as the telescoping series
\begin{equation}
\mathcal{L}^{(i)} = \sum_{k=0}^\infty \mathcal{L}^{(i)}_k - \mathcal{L}^{(i)}_{k-1} = \sum_{k=0}^\infty \Delta^{(i)}_k.
\end{equation}
Now we may write the tensor product of the exact operators as the tensor product of the telescoping sums, and interchange the product and sum:
\begin{eqnarray}
\mathcal{L}^{(1)} \otimes \cdots \otimes \mathcal{L}^{(d)} &= & \sum_{k_1=0}^\infty \Delta^{(1)}_{k_1} \otimes \cdots \otimes \sum_{k_d=0}^\infty \Delta^{(d)}_{k_d} \nonumber \\
 &= & \sum_{\mathbf{k} = 0}^\infty \Delta_{k_1}^{(1)} \otimes \cdots \otimes \Delta_{k_d}^{(d)}
\end{eqnarray}
\noindent Smolyak's idea is to approximate the tensor product operator with truncations of this sum \cite{Smolyak1963}:
\begin{equation}
\label{eq:generalSmolyakDiff}
A(\mathcal{K},d,\mathcal{L}) := \sum_{\mathbf{k} \in \mathcal{K}} \Delta_{k_1}^{(1)} \otimes \cdots \otimes \Delta_{k_d}^{(d)} .
\end{equation}
% 
% where $\vec{L}= \left(L^{(1)}, \ldots, L^{(d)}\right)$. 
We refer to the multi-index set $\mathcal{K}$ as the \emph{Smolyak multi-index set}, and it must be admissible for the sum to telescope correctly. Smolyak specifically suggested truncating with a total order multi-index set, which is the most widely studied choice. However, we can compute the approximation with any admissible multi-index set. Although the expression above is especially clean, it is not the most useful form for computation. We can reorganize the terms of (\ref{eq:generalSmolyakDiff}) to construct a weighted sum of the tensor operators:
\begin{equation}
\label{eq:generalSmolyakWeighted}
A(\mathcal{K}, d, \mathcal{L}) = \sum_{\mathbf{k} \in \mathcal{K}} c_\mathbf{k} \, \mathcal{L}^{(1)}_{k_1} \otimes \cdots \otimes \mathcal{L}^{(d)}_{k_d} ,
\end{equation}
where $c_\mathbf{k}$ are integer \emph{Smolyak coefficients} computed from the combinatorics of the difference formulation. One can compute the coefficients through a simple iteration over the index set and use (\ref{eq:generalSmolyakDiff}) to determine which full tensor rules are incremented or decremented. In general, these coefficients are non-zero near the leading surface of the Smolyak multi-index set, reflecting the mixing of the most accurate constituent full tensor approximations.

If each sequence of one-dimensional operators converges, then the Smolyak approximation converges 
{to the tensor product of exact operators} % could omit this line
as $\mathcal{K} \to \mathbb{N}^d_0$. For the isotropic simplex index set, some precise rates of convergence are known with respect to the side length of the simplex \cite{Wasilkowski1999,Wasilkowski2005a,Wasilkowski1995,Sickel2007a,Sickel2009}. Although general admissible Smolyak multi-index sets are difficult to study theoretically, they allow detailed customization to the anisotropy of a particular function. 

\subsection{Exactness of Smolyak algorithms}
In the one-dimensional and full tensor settings, we have characterized approximation algorithms through their exact sets---those inputs for which the algorithm is precise. This section shows that if the constituent one-dimensional approximations have nested exact sets, Smolyak algorithms are the ideal blending of different full tensor approximations from the perspective of exact sets; that is, the exact set of the Smolyak algorithm contains the union of the exact sets of the component full tensor approximations. This result will facilitate subsequent analysis of sparse quadrature and pseudospectral approximation algorithms. This theorem and our proof closely follow the framework provided by Novak and Ritter \cite{Novak1996,Novak1999a,Barthelmann2000}, but include a generalization to arbitrary Smolyak multi-index sets.

\begin{theorem}
\label{thm:SmolyakAccuracy}
Let $A(\mathcal{K}, d, \mathcal{L})$ be a Smolyak algorithm composed of linear operators with nested exact sets, i.e., with $m \leq m^\prime$ implying that $\mathcal{E} (\mathcal{L}^{(i)}_m ) \subseteq \mathcal{E} ( \mathcal{L}^{(i)}_{m^\prime})$ for $i=1 \ldots d$, where $\mathcal{K}$ is admissible. Then the exact set of $A(\mathcal{K}, d, \mathcal{L})$ contains
\begin{eqnarray}
\mathcal{E}\left (A(\mathcal{K}, d,\mathcal{L})\right ) &\supseteq &\bigcup_{\mathbf{k} \in \mathcal{K}} \mathcal{E} \left ( \mathcal{L}^{(1)}_{k_1} \otimes \cdots \otimes \mathcal{L}^{(d)}_{k_d} \right ) \nonumber \\
&\supseteq & \bigcup_{\mathbf{k} \in \mathcal{K}} \mathcal{E}(\mathcal{L}^{(1)}_{k_1}) \otimes \cdots \otimes \mathcal{E}(\mathcal{L}^{(d)}_{k_d}) .
\end{eqnarray}

\end{theorem}

\proof{
We begin by introducing notation to incrementally build a multi-index set dimension by dimension. For a multi-index set $\mathcal{K}$ of dimension $d$, let the restriction of the multi-indices to the first $i$ dimensions be $\mathcal{K}^{(i)} := \{\mathbf{k}_{1:i} = (k_1, \ldots, k_i) : \mathbf{k} \in \mathcal{K}\}$. Furthermore, define subsets of $\mathcal{K}$ based on the $i$\textsuperscript{th} element of the multi-indices, $\mathcal{K}^{(i)}_j := \{\mathbf{k}_{1:i} : \mathbf{k} \in \mathcal{K}^{(i)} \ \mathit{and} \ {k}_{i+1} = j\}$. These sets are nested, $\mathcal{K}^{(i)}_j \supseteq \mathcal{K}^{(i)}_{j+1}$, because $\mathcal{K}$ is admissible. 
Also let ${k}^\mathrm{max}_i$ denote the maximum value of the $i$\textsuperscript{th} component of the multi-indices in the set $\mathcal{K}$.

Using this notation, one can construct $\mathcal{K}$ inductively, 
\begin{eqnarray}
\mathcal{K}^{(1)} &= &  \{1, \ldots,  {k}^\mathrm{max}_{1}\}\\
\mathcal{K}^{(i)} &= & \bigcup_{j = 1}^{{k}^\mathrm{max}_{i}} \mathcal{K}^{(i-1)}_j \otimes j, \ \ i = 2 \ldots d. \label{e:inductK}
\end{eqnarray}

It is sufficient to prove that the Smolyak operator is exact for an arbitrary $f$ with tensor structure, $f = f_1 \times \cdots \times f_{d}$. Suppose there exists a $\mathbf{k}^\ast$ such that  $f \in \mathcal{E}(\mathcal{L}^{(\mathbf{d})}_{\mathbf{k}^\ast})$. We will show that if $\mathcal{K}$ is an admissible multi-index set containing $\mathbf{k}^\ast$, then $A (\mathcal{K}, d,\mathcal{L})$ is exact on $f$. We do so by induction on the dimension $i$ of the Smolyak operator and the function. 

First, consider the $i=1$ case. $A(\mathcal{K}^{(1)}, 1,\mathcal{L}) = \mathcal{L}^{(1)}_{{k}^\mathrm{max}_1}$, where ${k}^\mathrm{max}_1 \geq {k}^\ast_1$. Hence $\mathcal{E}(A(\mathcal{K}^{(1)}, 1,\mathcal{L})) = \mathcal{E}(\mathcal{L}^{(1)}_{{k}^\mathrm{max}_1}).$

For the induction step, we construct the $(i+1)$-dimensional Smolyak operator in terms of the $i$-dimensional operator: 
\begin{equation}
A(\mathcal{K}^{(i+1)}, i+1, \mathcal{L}) = \sum_{j = 1}^{{k}^\mathrm{max}_{i+1}} A(\mathcal{K}^{(i)}_j, i, \mathcal{L}) \otimes (\mathcal{L}^{(i+1)}_{j}-\mathcal{L}^{(i+1)}_{j-1}) .
\label{e:fullsum}
\end{equation}
This sum is over increasing levels of accuracy in the $i+1$ dimension. We know the level required for the approximate operator to be exact  in this dimension; this may be expressed as
\begin{equation}
\mathcal{L}^{(i+1)}_{j}(f_{i+1}) = \mathcal{L}^{(i+1)}_{j-1}(f_{i+1}) = \mathcal{L}^{(i+1)}(f_{i+1}) \ \mathrm{when} \ j-1 \geq {k}^\ast_{i+1} .
\end{equation}
Therefore the sum (\ref{e:fullsum}) can be truncated at the ${k}^\ast_{i+1}$ term, as the differences of higher terms are zero when applied to $f$:
\begin{equation}
A(\mathcal{K}^{(i+1)}, i+1, \mathcal{L}) = \sum_{j = 1}^{{k}^\ast_{i+1}} A(\mathcal{K}^{(i)}_j, i,\mathcal{L}) \otimes (\mathcal{L}^{(i+1)}_{j}-\mathcal{L}^{(i+1)}_{j-1}).
\end{equation}
Naturally, $\mathbf{k}^\ast_{1:i} \in \mathcal{K}^{(i)}_{{k}^\ast_{i+1}}$. By nestedness,  $\mathbf{k}^\ast_{1:i}$ is also contained in $\mathcal{K}^{(i)}_{j}$ for $j \leq {k}^\ast_{i+1}$. The induction hypothesis then guarantees
\begin{equation}
f_1 \otimes \cdots \otimes f_i \in \mathcal{E}(A(\mathcal{K}^{(i)}_j, i, \mathcal{L})), \ \forall j \leq {k}^\ast_{i+1} .
\end{equation}
Applying the $(i+1)$-dimensional Smolyak operator to the truncated version of $f$ yields
\begin{eqnarray}
& & A(\mathcal{K}^{(i+1)}, i+1, \mathcal{L})(f_1 \otimes \cdots \otimes f_{i+1}) \nonumber \\
& = &  \sum_{j = 1}^{{k}^\ast_{i+1}} A(\mathcal{K}^{(i)}_j, i, \mathcal{L})(f_1 \otimes \cdots \otimes f_i) \otimes (\mathcal{L}^{(i+1)}_{j}-\mathcal{L}^{(i+1)}_{j-1})(f_{i+1}) .
\end{eqnarray}
Since each of the $i$-dimensional Smolyak algorithms is exact, by the induction hypothesis, we replace them with the true operators and rearrange by linearity to obtain
\begin{eqnarray}
A(\mathcal{K}^{(i+1)}, i+1, \mathcal{L})(f_1 \otimes \cdots \otimes f_{i+1})  &= & \mathcal{L}^{(\mathbf{i})}(f_1 \otimes \cdots \otimes f_i)\otimes \sum_{j = 1}^{{k}^\ast_{i+1}} (\mathcal{L}^{(i+1)}_{j}-\mathcal{L}^{(i+1)}_{j-1})(f_{i+1}) \nonumber \\
&= & \mathcal{L}^{(\mathbf{i})}(f_1 \otimes \cdots \otimes f_i)\otimes \mathcal{L}^{(i+1)}_{{k}^\ast_{i+1}}(f_{i+1}). \label{e:lastdim}
\end{eqnarray}
The approximation in the $i+1$ dimension is exactly of the level needed to be exact on the $(i+1)$\textsuperscript{th} component of $f$. Then (\ref{e:lastdim}) becomes
\begin{equation}
\mathcal{L}^{(\mathbf{i})}(f_1 \otimes \cdots \otimes f_i)\otimes \mathcal{L}^{(i+1)}(f_{i+1}) = \mathcal{L}^{(\mathbf{i+1})}(f_1 \otimes \cdots \otimes f_{i+1})
\end{equation}
Thus the Smolyak operator is precise for $f$, and the claim is proven.} \endproof %end proof

\subsection{Smolyak quadrature}
We recall the most familiar use of Smolyak algorithms, sparse quadrature. Consider a family of one-dimensional quadrature rules $ \mathcal{Q}_k^{(i)} $ in each dimension $i=1 \ldots d$; denote these rules by ${\mathcal{Q}}$. The resulting Smolyak algorithm is written as: 
\begin{equation}
A(\mathcal{K}, d, {\mathcal{Q}}) = \sum_{\mathbf{k} \in \mathcal{K}} c_\mathbf{k} \mathcal{Q}^{(\mathbf{d})}_\mathbf{k}.
\end{equation}

\begin{figure}
\centering

\subfloat[The exact set for a level-four Smolyak quadrature in two dimensions, based on \emph{linear} growth Gaussian quadrature rules.]
{
\label{fig:smolyakQuadratureAccuraciesA}
\includegraphics[scale=.7]{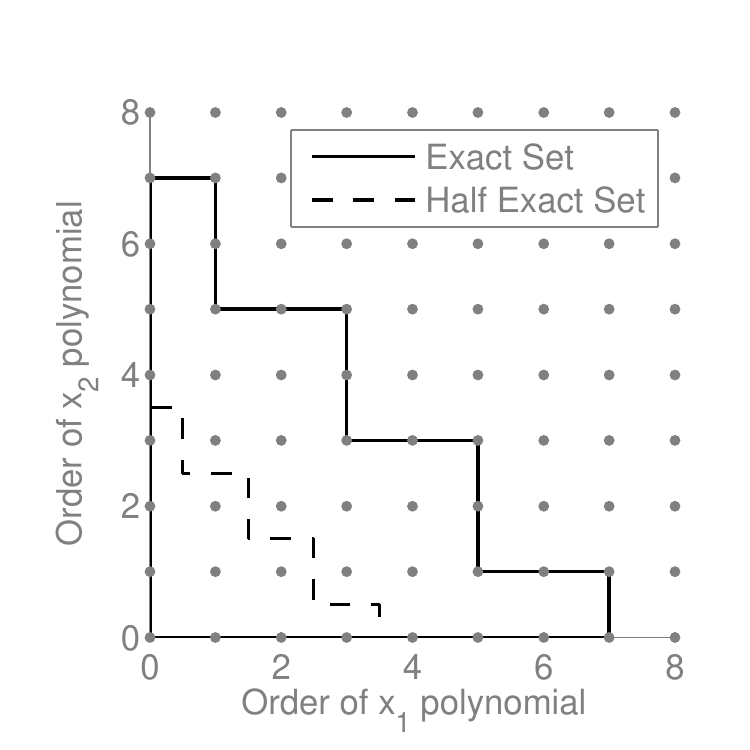}
}
\qquad
\subfloat[The exact set for a level-three Smolyak quadrature in two dimensions, based on \emph{exponential} growth Gaussian quadrature rules.]
{
\label{fig:smolyakQuadratureAccuraciesB}
\includegraphics[scale=.7]{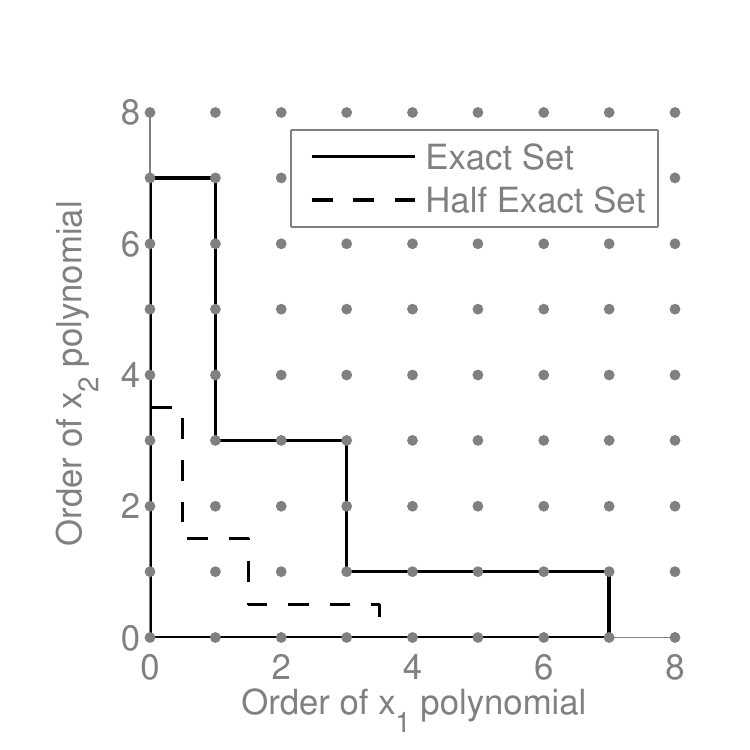}
}
\caption{The exact set diagram for two Smolyak quadrature rules, and the corresponding basis for a Smolyak pseudospectral approximation. $\mathcal{E}(\mathcal{Q})$ is shown in a solid line, $\mathcal{E}_2(\mathcal{Q})$ is the dashed line. The staircase appearance results from the superposition of rectangular full tensor exact sets.}
\label{fig:smolyakQuadratureAccuracies}
\end{figure}

This approximation inherits its convergence from the one-dimensional operators. The set of functions that are exactly integrated by a Smolyak quadrature algorithm is described as a corollary of Theorem \ref{thm:SmolyakAccuracy}.
\begin{corollary}
\label{thm:sparseQuadAccuracy}
For a sparse quadrature rule satisfying the hypotheses of Theorem \ref{thm:SmolyakAccuracy},
\begin{equation}
\mathcal{E} \left (A (\mathcal{K},d,{\mathcal{Q}}) \right ) \supseteq 
\bigcup_{\mathbf{k} \in \mathcal{K}} \mathcal{E}(\mathcal{Q}^{(\mathbf{d})}_\mathbf{k})  
\end{equation}
\end{corollary}
Quadrature rules with polynomial accuracy do have nested exact sets, as required by the theorem. An example of Smolyak quadrature exact sets is shown in Figure \ref{fig:smolyakQuadratureAccuracies}. 

%Additionally, other separable functions may be precisely integrated.
%%
%\begin{Lemma}
%\label{thm:sparseSeparable}
%If $f$ is separable such that  $f = f_1(x_1)f'(x_2, \ldots, x_d)$ and $\mathcal{I}^{(1)}(f_1) = 0$, and if for every full tensor rule $\mathbf{k} \in \mathcal{K}$ we have $f_1 \in \mathcal{E}(\mathcal{Q}^{(1)}_{{k}_1})$, then $f \in \mathcal{E}(A(\mathcal{K},d,{\mathcal{Q}}))$.
%\end{Lemma}
%
%\proof{
%If every full tensor rule in the Smolyak algorithm correctly computes that $f_1$ integrates to zero, then their weighted sum is necessarily zero.
%}\endproof

\subsection{Smolyak pseudospectral approximation}
Applying Smolyak's algorithm to pseudospectral approximation operators yields a sparse algorithm that converges under similar conditions as the one-dimensional operators from which it is constructed. This algorithm is written as
\begin{equation}
A(\mathcal{K}, d, {\mathcal{S}}) = \sum_{\mathbf{k} \in \mathcal{K}} c_\mathbf{k} \mathcal{S}^{(\mathbf{d})}_\mathbf{k} .
\end{equation}
The Smolyak algorithm is therefore a sum of different full tensor pseudospectral approximations, where each approximation is built around the polynomial accuracy of a single full tensor quadrature rule. It is not naturally expressed as a set of formulas for the polynomial coefficients, because different approximations include different polynomials. The term $\Psi_\mathbf{j}$ is included in the Smolyak approximation if and only if $\exists \mathbf{k} \in \mathcal{K}: \Psi_\mathbf{j} \in \mathcal{E}_2(\mathcal{Q}^{(\mathbf{d})}_\mathbf{k})$. Here, $\mathcal{Q}^{(\mathbf{d})}_\mathbf{k}$ is the full tensor quadrature rule used by the full tensor pseudospectral approximation $\mathcal{S}^{(\mathbf{d})}_\mathbf{k}$. As in the full tensor case, the half exact set of a Smolyak quadrature rule defines the range of the Smolyak pseudospectral approximation.

Once again, the Smolyak construction guarantees that the convergence of this approximation is inherited from its constituent one-dimensional approximations. Our choices for the pseudospectral operators ensure nestedness of the constituent exact sets, so we may use Theorem~\ref{thm:SmolyakAccuracy} to ensure that Smolyak pseudospectral algorithms are exact on their range. 

\begin{corollary}
If the constituent one-dimensional pseudospectral rules have no internal aliasing and satisfy the conditions of Theorem \ref{thm:SmolyakAccuracy}, then the resulting Smolyak pseudospectral algorithm has no internal aliasing.
\end{corollary}

We additionally provide a theorem that characterizes the external aliasing properties of Smolyak pseudospectral approximation, which the next section will contrast with direct quadrature. 

\begin{theorem}
\label{thm:smolyakExternal}
Let $\Psi_\mathbf{j}$ be a polynomial term included in the expansion provided by the Smolyak algorithm $A(\mathcal{K}, d, {\mathcal{S}})$, and let $\Psi_{\mathbf{j}^{\prime}}$ be a polynomial term not included in the expansion. There is no external aliasing of $\Psi_{\mathbf{j}^{\prime}}$ onto $\Psi_\mathbf{j}$ if any of the following conditions is satisfied: (a) there exists a dimension $i$ for which ${j}^{\prime}_i < {j}_i$; or (b) there exists a multi-index $\mathbf{k} \in \mathcal{K}$ such that $\Psi_\mathbf{j}$ is included in the range of  $\mathcal{S}^{(\mathbf{d})}_\mathbf{k}$ and $\Psi_{\mathbf{j}^{\prime}} \Psi_\mathbf{j} \in \mathcal{E}(\mathcal{Q}^{(\mathbf{d})}_\mathbf{k})$, where $\mathcal{Q}^{(\mathbf{d})}_\mathbf{k}$ is the quadrature rule used in $\mathcal{S}^{(\mathbf{d})}_\mathbf{k}$.
\end{theorem}

\proof{If condition (a) is satisfied, then $\Psi_\mathbf{j}$ and $\Psi_{\mathbf{j}^{\prime}}$ are orthogonal in dimension $i$, and hence that inner product is zero. Every quadrature rule that computes the coefficient $f_\mathbf{j}$ corresponding to basis term $\Psi_\mathbf{j}$ is accurate for polynomials of at least order $2\mathbf{j}$. Since ${j}^{\prime}_i +  {j}_i < 2{j}_i$, every rule that computes the coefficient can numerically resolve the orthogonality, and therefore there is no aliasing. If condition (b) is satisfied, then the result follows from the cancellations exploited by the Smolyak algorithm, as seen in the proof of Theorem \ref{thm:SmolyakAccuracy}.} \endproof

These two statements yield extremely useful properties. First, any Smolyak pseudospectral algorithm, regardless of the admissible Smolyak multi-index set used, has no internal aliasing; this feature is important in practice and not obviously true. Second, while there is external aliasing as expected, the algorithm uses basis orthogonality to limit which external coefficients can alias onto an included coefficient. The Smolyak pseudospectral algorithm is thus a practically ``useful'' approximation, in that one can tailor it to perform a desired amount of work while guaranteeing reliable approximations of the selected coefficients. Computing an accurate approximation of the function only requires including sufficient terms so that the truncation and external aliasing errors are small.

\section{Comparing direct quadrature to Smolyak pseudospectral approximation}
\label{sec:comparison}
The current UQ literature often suggests a \emph{direct quadrature} approach for constructing polynomial chaos expansions \cite{Xiu2009,LeMaitre2010,Eldred2009:local,Huan2012}. In this section, we describe this approach and show that, in comparison to a true Smolyak algorithm, it is inaccurate or inefficient in almost all cases. Our comparisons will contrast the theoretical error performance of the algorithms and provide simple numerical examples that illustrate typical errors and why they arise.

\subsection{Direct quadrature polynomial expansions}

At first glance, direct quadrature is quite simple. First, choose a multi-index set $\mathcal{J}$ to define a truncated polynomial expansion:
\begin{equation}
f \approx \sum_{\mathbf{j \in \mathcal{J}}} \tilde{f}_\mathbf{j} \Psi_\mathbf{j} .
\end{equation}
The index set $\mathcal{J}$ is typically admissible, but need not be. Second, select any $d$-dimensional quadrature rule $\mathcal{Q}^{(\mathbf{d})}$, and estimate every coefficient as:
\begin{equation}
\tilde{f_\mathbf{j}}= \mathcal{Q}^{(\mathbf{d})}(f\Psi_\mathbf{j}) .
\end{equation}
Unlike the Smolyak approach, we are left to choose $\mathcal{J}$ and $\mathcal{Q}^{(\mathbf{d})}$ independently, giving the appearance of flexibility. In practice, this produces a more complex and far more subtle version of the truncation trade-off discussed in Section \ref{sec:approximations}. Below, we will be interested in selecting $\mathcal{Q}$ and $\mathcal{J}$ to replicate the quadrature points and output range of the Smolyak approach, as it provides a benchmark for achievable performance.

Direct quadrature does not converge for every choice of $\mathcal{J}$ and $\mathcal{Q}^{(\mathbf{d})}$; consider the trivial case where $\mathcal{J}$ does not grow infinitely. It is possible that including far too many terms in $\mathcal{J}$ relative to the polynomial exactness of $\mathcal{Q}^{(\mathbf{d})}$ could lead to a non-convergent algorithm. Although this behavior contrasts with the straightforward convergence properties of Smolyak algorithms, most reasonable choices for direct quadrature do converge, and hence this is not our primary argument against the approach. 

Instead, our primary concern is aliasing in direct quadrature and how it reduces performance at finite order. Both internal and external aliasing are governed by the same basic rule, which is just a restatement of how we defined aliasing in Section \ref{sec:aliasing}.

\begin{Remark}
\label{rem:dqAliasing}
For a multi-index set $\mathcal{J}$ and a quadrature rule $\mathcal{Q}^{(\mathbf{d})}$, the corresponding direct quadrature polynomial expansion has no aliasing between two polynomial terms if $\Psi_\mathbf{j} \Psi_\mathbf{j^\prime} \in \mathcal{E}(\mathcal{Q}^{(\mathbf{d})})$.
\end{Remark}

The next two sections compare the internal and external aliasing with both theory and simple numeric examples.

\subsection{Internal aliasing in direct quadrature}

As an extension of Remark \ref{rem:dqAliasing}, direct quadrature has no internal aliasing whenever every pair $\mathbf{j},\mathbf{j^\prime} \in \mathcal{J}$ has no aliasing. We can immediately conclude that for any basis set $\mathcal{J}$, there is some quadrature rule sufficiently powerful to avoid internal aliasing errors. In practice, however, this rule may not be a desirable one.

\begin{Example}
Assume that for some function with two inputs, we wish to include the polynomial basis terms $(a,0)$ and $(0,b)$. By Remark \ref{rem:dqAliasing}, the product of these two terms must be in the exact set; hence, the quadrature must include at least a full tensor rule of accuracy $(a,b)$. Although we have not asked for any coupling, direct quadrature must assume full coupling of the problem in order to avoid internal aliasing.
\end{Example}

Therefore we reach the surprising conclusion that direct quadrature inserts significant coupling into the problem, whereas we selected a Smolyak quadrature rule in hopes of leveraging the absence of that very coupling---making the choice inconsistent. For most sparse quadrature rules, we cannot include as many polynomial terms as the Smolyak pseudospectral approach without incurring internal aliasing, because the quadrature is not powerful enough in the mixed dimensions. 

\begin{figure}[htb]
	
	\centering
	\includegraphics[scale=.5]{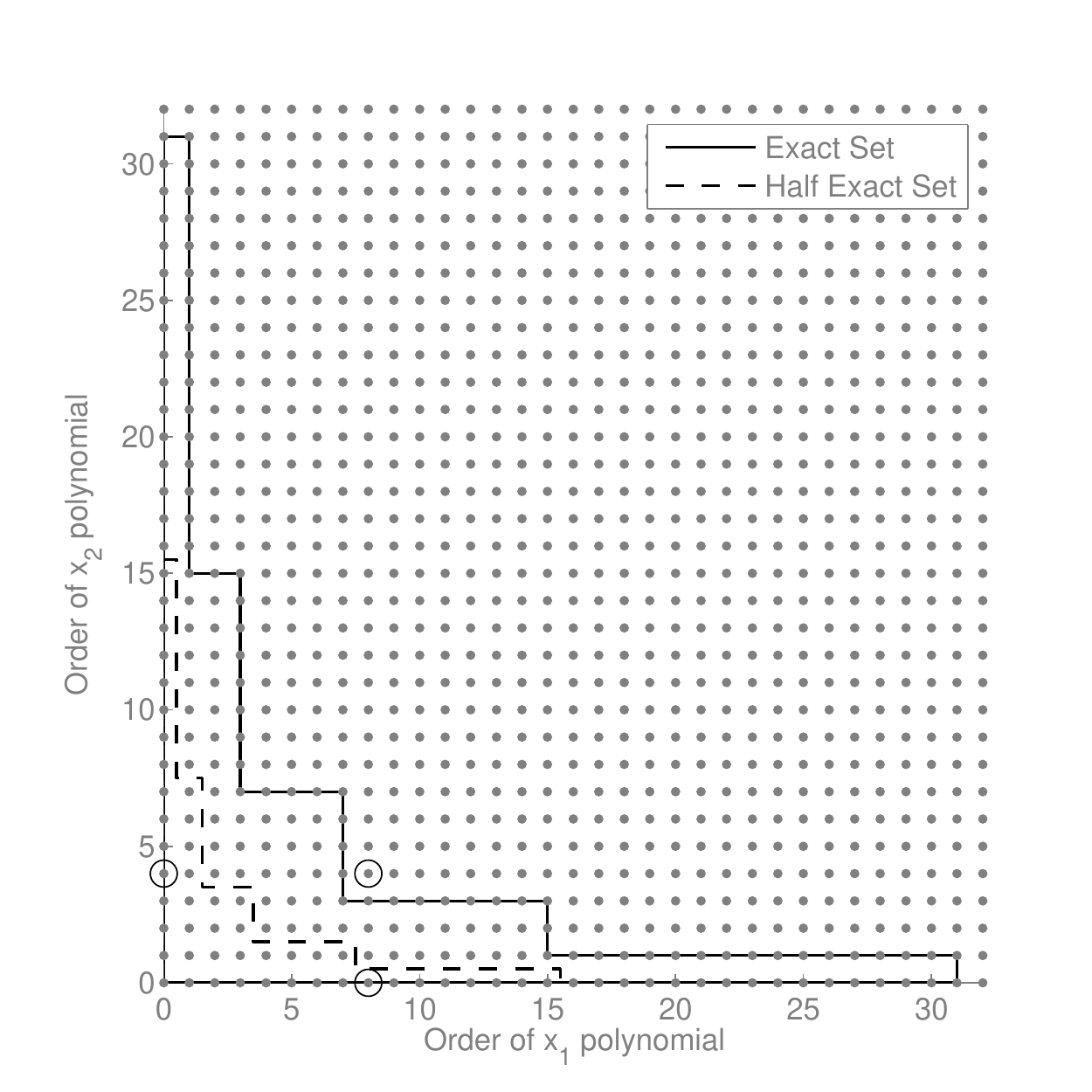}
	\caption{The exact set and polynomials included in the direct quadrature construction from Example \ref{ex:directAliasing}.}
	\label{fig:ExpDirectQuadrature}
\end{figure}

\begin{figure}
\centering

\subfloat[Smolyak Pseudospectral]
{
\includegraphics[scale=.7]{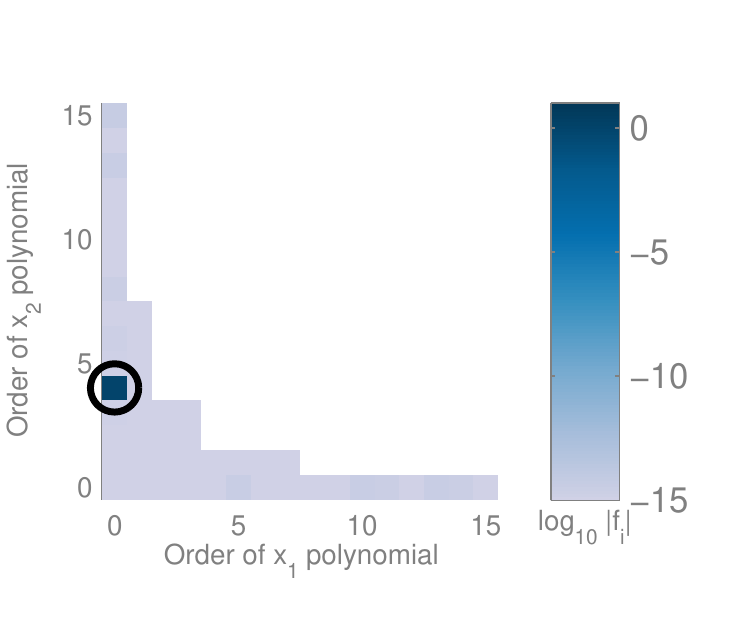}
}
\qquad
\subfloat[Direct Quadrature]
{
\includegraphics[scale=.7]{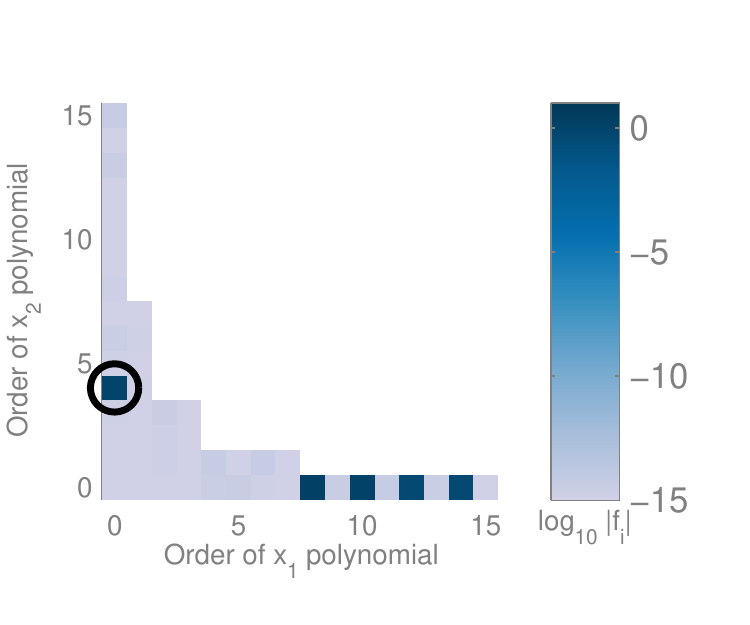}
}
\caption{Numerical results for Example \ref{ex:directAliasing}; each color square indicates the log of the coefficient magnitude for the basis function at that position. The circle identifies the correct non-zero coefficient.
}
\label{fig:ExpInternalAliasing}
\end{figure}

\begin{Example}
\label{ex:directAliasing}
Let $\mathbf{X}$ be the two-dimensional domain $[-1,1]^2$. Select a uniform weight function, which corresponds to a Legendre polynomial basis. Let $f(x,y) = \psi_0(x)\psi_4(y)$. Use Gauss-Legendre quadrature and an exponential growth rule, such that $p^{(i)}(m) = 2^{m-1}$. Select a sparse quadrature rule based on a total order multi-index set $\mathcal{Q}^2_{\mathcal{K}^{t}_5}$. Figure \ref{fig:ExpDirectQuadrature} shows the exact set of this Smolyak quadrature rule (solid line) along with its half-exact set (dashed line), which encompasses all the terms in the direct quadrature polynomial expansion.

Now consider the $\mathbf{j} = (8,0)$ polynomial, which is in the half-exact set. The product of the $(0,4)$ and $(8,0)$ polynomial terms is $(8,4)$, which is not within the exact set of the sparse rule. Hence, $(0,4)$ aliases onto $(8,0)$ because this quadrature rule has limited accuracy in the mixed dimensions. 

Using both the Smolyak pseudospectral and direct quadrature methods, we numerically compute the polynomial expansion for this example. The resulting coefficients are shown in Figure \ref{fig:ExpInternalAliasing}. Even though the two methods use the same information and project $f$ onto the same basis, the Smolyak result has no internal aliasing while direct quadrature shows significant internal aliasing. Although both methods correctly compute the $(0,4)$ coefficient, direct quadrature shows aliasing on $(8,0)$ as predicted, and also on $(10,0)$, $(12,0)$, and $(14,0)$. In this case, direct quadrature is unable to determine the order of the input function or even whether the input is function of $x_1$ or $x_2$. Alternating terms are computed correctly because of the parity of the functions.
\end{Example}

The \OO{1} errors observed in this simple example demonstrate why it is crucial to eliminate internal aliasing in the construction of one-dimensional pseudospectral approximations and to ensure that the full tensor and Smolyak algorithms inherit that property. More complex functions demonstrate the same type of error, except that the errors resulting from multiple source terms are superimposed. Examples of the latter are given by Constantine \textit{et al.} \cite{Constantine2012}.  

There are some choices for which direct quadrature has no internal aliasing: full tensor quadrature rules and, notably, Smolyak quadrature constructed from one-dimensional Gaussian quadrature rules with $p^{(i)}(m) = m$, truncated according to an isotropic total-order multi-index set. However, many useful sparser or more tailored Smolyak quadrature rules, e.g., based on exponential growth quadrature rules or adaptive anisotropic Smolyak index sets, will incur internal aliasing if the basis selection matches the range of the Smolyak algorithm. This makes them a poor choice when a comparable Smolyak pseudospectral algorithm uses the same evaluation points and produces an approximation with the same polynomial terms, but is guaranteed by construction to have zero internal aliasing. Alternately, it is possible to select a sufficiently small polynomial basis to avoid internal aliasing, but this approach requires unnecessary conservatism that could easily be avoided with a Smolyak pseudospectral approximation.

\subsection{External aliasing}
The difference in external aliasing between direct quadrature and Smolyak pseudospectral approximation is much less severe. Both methods exhibit external aliasing from terms far outside the range of the approximation, as such errors are a necessary consequence of using finite order quadrature. Since the methods are constructed from similar constituent one-dimensional quadrature rules, aliasing is of similar magnitude when it occurs. 

Comparing Theorem \ref{thm:smolyakExternal}, condition \textit{(b}), and Remark \ref{rem:dqAliasing}, we observe that if the direct quadrature method has no external aliasing between two basis terms, the equivalent Smolyak pseudospectral algorithm will not either. Yet the two methods perform differently because of their behavior on separable functions. Condition \textit{(a)} of Theorem \ref{thm:smolyakExternal} provides an additional condition under which external aliasing will not occur under a Smolyak pseudospectral algorithm, and thus it has strictly less external aliasing in general.

\begin{figure}
\centering

\subfloat[Smolyak Pseudospectral]
{
\includegraphics[scale=.7]{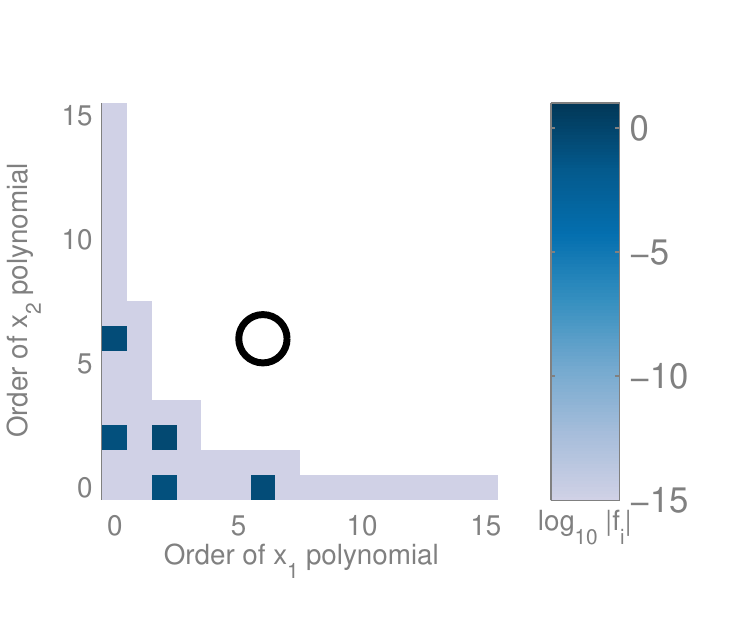}
}
\qquad
\subfloat[Direct Quadrature]
{
\includegraphics[scale=.7]{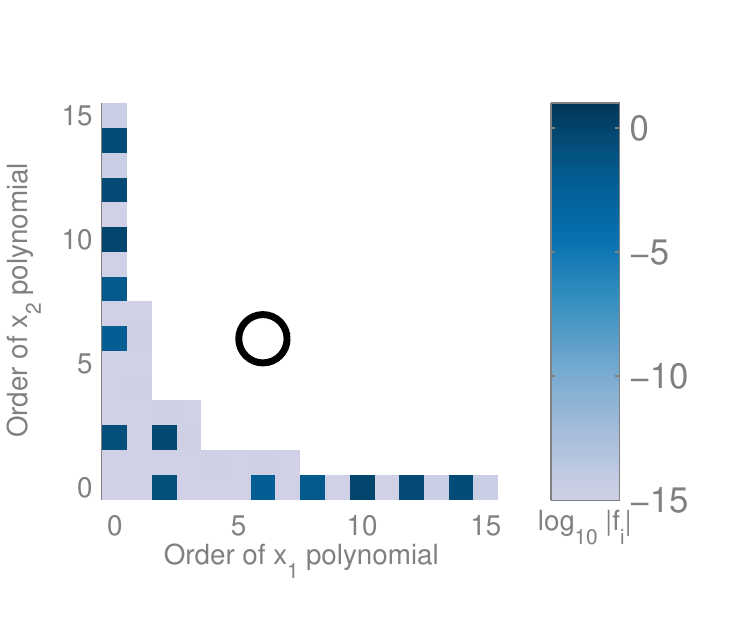}
}
\caption{Numerical results for Example \ref{ex:externalAliasing}; each color square indicates the log of the coefficient magnitude for the basis function at its position. The circle indicates the correct non-zero coefficient. The Smolyak pseudospectral approach has fewer terms corrupted by external aliasing in this case.
}
\label{fig:externalAliasing}
\end{figure}

\begin{Example}
\label{ex:externalAliasing}
If we repeat Example \ref{ex:directAliasing} but choose $f$ to be a polynomial outside the approximation space, $f = \psi_6(x)\psi_6(y)$, we obtain the results in Figure \ref{fig:externalAliasing}. Now every non-zero coefficient is the result of external aliasing. Direct quadrature correctly computes some terms because of either parity or the few cases where Remark \ref{rem:dqAliasing} is satisfied. However, the Smolyak approach has fewer errors because the terms not between $(0,0)$ and $(6,6)$ are governed by condition \textit{(a)} of Theorem \ref{thm:smolyakExternal}, and hence have no external aliasing. 
\end{Example}

This example is representative of the general case. Direct quadrature always incurs at least as much external aliasing as the Smolyak approach, and the methods become equivalent if the external term causing aliasing is of very high order. Although both methods will always exhibit external aliasing onto coefficients of the approximation for non-polynomial inputs, the truncation can in principle be chosen to include all the important terms, so that the remaining external aliasing is acceptably small.

\subsection{Summary of comparison}
% epitome?

Compared to the Smolyak pseudospectral approach, direct quadrature yields larger internal \textit{and} external aliasing errors. Because of these aliasing errors, direct quadrature is essentially unable to make efficient use of most sparse quadrature rules. The Smolyak pseudospectral approach, on the other hand, is guaranteed never to have internal aliasing if the one-dimensional pseudospectral operators are chosen according to simple guidelines. We therefore recommend against using direct quadrature. The remainder of the paper will focus on extensions of the basic Smolyak pseudospectral approach.

\section{Adaptive polynomial approximations}
\label{sec:adaptive}

When constructing a polynomial approximation of a black-box computational model, there are two essential questions: first, which basis terms should be included in the expansion; and second, what are the coefficients of those basis terms? The Smolyak construction allows detailed control over the truncation of the polynomial expansion and the work required to compute it. Since we typically do not have \textit{a priori} information about the functional relationship generated by a black-box model, we develop an adaptive approach to tailor the Smolyak approximation to this function, following the dimension-adaptive quadrature approaches of Gerstner \& Griebel~\cite{Gerstner2003} and Hegland~\cite{Hegland2003}.

The Smolyak algorithm is well suited to an adaptive approach. The telescoping sum converges in the same sense as the constituent one-dimensional operators as the index set grows to include $\mathbb{N}^d_0$, so we can simply add more terms to improve the approximation until we are satisfied. We separate our adaptive approach into two components: local improvements to the Smolyak multi-index set and a global stopping criterion.

\subsection{Dimension adaptivity}
Dimension adaptivity is responsible for identifying multi-indices to add to a Smolyak multi-index set in order to improve its accuracy. A standard refinement approach is simply to grow an isotropic simplex of side length $n$. \cite{Gerstner2003} and \cite{Hegland2003} instead suggest a series of greedy refinements that customize the Smolyak algorithm to a particular problem.

The refinement used in \cite{Gerstner2003} is to select a multi-index $\mathbf{k} \in \mathcal{K}$ and to add the forward neighbors of $\mathbf{k}$ that are admissible. The multi-index $\mathbf{k}$ is selected via an error indicator $\epsilon(\mathbf{k})$. We follow \cite{Gerstner2003} and assume that whenever $\mathbf{k}$ contributes strongly to the result of the algorithm, it represents a subspace that likely needs further refinement. 

Let $\mathbf{k}$ be a multi-index such that $\mathcal{K}^\prime := \mathcal{K} \cup \mathbf{k}$, where $\mathcal{K}$ and $\mathcal{K}^\prime$ are admissible multi-index sets. The triangle inequality (for some appropriate norm, see Section \ref{s:adaptivecomments}) bounds the change in the Smolyak approximation produced by adding $\mathbf{k}$ to $\mathcal{K}$, yielding a useful error indicator:
\begin{equation}
\label{eq:localError}
\|A(\mathcal{K}^\prime,d,\mathcal{L}) - A(\mathcal{K},d,\mathcal{L}) \| \leq \| \Delta_{{k}_1}^1 \otimes \cdots \otimes \Delta_{{k}_d}^d   \| =: \epsilon(\mathbf{k}) .
\end{equation}
Conveniently, this error indicator does not change as $\mathcal{K}$ evolves, so we need only compute it once. At each adaptation step, we find the $\mathbf{k}$ that maximizes $\epsilon(\mathbf{k})$ and that has at least one admissible forward neighbor. Then we add those forward neighbors. 

\subsection{Termination criterion}

Now that we have a strategy to locally improve a multi-index set, it is useful to have a global estimate of the error of the approximation, $\epsilon_g$. We cannot expect to compute the exact error, but even a rough estimate is useful. We follow Gerstner \& Griebel's choice of global error indicator
\begin{equation}
\label{eq:globalError}
\epsilon_g := \sum \epsilon(\mathbf{k}) ,
\end{equation}
where the sum is taken over all the multi-indices that are eligible for local adaptation at any particular step (i.e., that have admissible forward neighbors) \cite{Gerstner2003}. The algorithm may be terminated when a particular threshold of the global indicator is reached, or when it falls by a specified amount. 

\subsection{Error indicators and work-considering algorithms}
\label{s:adaptivecomments}

Thus far we have presented the adaptation strategy without reference to the problem of polynomial approximation. In this specific context, we use the $L^2(\mathbf{X}, w)$ norm in (\ref{eq:localError}), because it corresponds to the convergence properties of pseudospectral approximation and thus seems an appropriate target for greedy refinements. This choice is a heuristic to accelerate performance---albeit one that is simple and natural, and has enjoyed success in numerical experiments (see Section~\ref{sec:experiments}). Moreover, the analysis of external aliasing in Theorem \ref{thm:smolyakExternal} suggests that, in the case of pseudospectral approximation, significant missing polynomial terms alias onto some of the \textit{included} lower-order coefficients, giving the algorithm a useful indication of which direction to refine. This behavior helps reduce the need for smoothness in the coefficient pattern. Section \ref{sec:quadRules} provides a small fix that further helps with even or odd functions.

One is not required to use this norm to define $\epsilon(\mathbf{k})$, however, and it is possible that other choices could serve as better heuristics for some problems. Unfortunately, making definitive statements about the properties or general utility of these heuristic refinement schemes is challenging. The approach described above is intended to be broadly useful, but specific applications may require experimentation to find better choices. 

% Though seemingly reasonable, these are not unique or optimal choices, and there certainly exist examples for which these strategies reduce performance. 
Beyond the choice of norm, a commonly considered modification to $\epsilon (\mathbf{k})$ is to incorporate a notion of the computational effort required to refine $\mathbf{k}$. Define $n(\mathbf{k})$ as the amount of work to refine the admissible forward neighbors of $\mathbf{k}$, e.g., the number of new function evaluation points. \cite{Gerstner2003} discusses an error indicator that provides a parameterized sliding scale between selecting the term with highest $\epsilon (\mathbf{k})$ and the lowest $n(\mathbf{k})$:
\begin{equation}
\epsilon_{w,1}(\mathbf{k}) = \max \left\{ w \frac{\epsilon (\mathbf{k})}{\epsilon (\mathbf{1})}, (1 - w)\frac{n (\mathbf{1})}{n (\mathbf{k})}\right\}.
\end{equation}
Here $w \in [0,1]$ is the tuning parameter, and $\epsilon (\mathbf{1})$ and $n (\mathbf{1})$ are the indicator and cost of the first term. Putting $w=0$ considers only the standard error indicator and $w=1$ considers only the cost. A different indicator with a similar intent is
\begin{equation}
\epsilon_{w,2}(\mathbf{k}) = \epsilon (\mathbf{k}) - \tilde{w} n(\mathbf{k}),
\end{equation}
where $\tilde{w}>0$ describes a conversion between error and work. Both of these methods will sometimes select terms of low cost even if they do not appear to provide immediate benefit to the approximation. Yet we find both methods to be challenging to use in practice, because of the difficulty of selecting the tuning parameter. One can remove this particular tuning requirement by taking a ratio:
\begin{equation}
\epsilon_{w,3}(\mathbf{k}) = \frac{\epsilon (\mathbf{k}) }{n(\mathbf{k})}.
\end{equation}
This indicator looks for ``efficient'' terms to refine---ones that are expected to yield greater error reduction at less cost---rather than simply the highest-error directions. We performed some numerical experiments with these methods, but none of the examples demonstrated significant improvement. Furthermore, poor choices of tuning parameters can be harmful because they can essentially make the algorithm revert to a non-adaptive form. We do not give detailed results here for those experiments because they are not particularly conclusive; some types of coefficient patterns may benefit from work-considering approaches, but this remains an open problem.

On a similar note, using $\epsilon_g$ as a termination criteria is also a heuristic. As our experiments in Section \ref{sec:experiments} will show, for most smooth functions $\epsilon_g$ is an excellent estimate of the approximation accuracy. In other cases, the indicator can be quite poor; hence one should not rely on it exclusively. In practice, we typically terminate the algorithm based on a combination of elapsed wall clock time, the global error indicator, and an error estimate computed from limited \textit{ad hoc} sampling.

\section{Numerical experiments}
\label{sec:experiments}
Our numerical experiments focus on evaluating the performance of
different quadrature rules embedded within the Smolyak pseudospectral
scheme, and on evaluating performance of the adaptive Smolyak
approximation strategy. Aside from the numerical examples of
Section~\ref{sec:comparison}, we do not investigate the performance of
direct quadrature any further. Given our theoretical analysis of
aliasing errors and the numerical demonstrations in
\cite{Constantine2012}, one can conclude without further demonstration
that destructive internal aliasing indeed appears in practice.

This section begins by discussing practical considerations in the
selection of quadrature rules. Then we evaluate convergence of Smolyak
pseudospectral approximation schemes (non-adaptive and adaptive) on
the Genz test functions. Next, we approximate a larger chemical
kinetic system, illustrating the efficiency and accuracy of the
adaptive method. Finally, we evaluate the quality of the global error
indicator on all of these examples.

\subsection{Selection of quadrature rules}
\label{sec:quadRules}
Thus far we have sidestepped practical questions about which quadrature rules exist or are most efficient. Our analysis has relied only on polynomial accuracy of quadrature rules; all quadrature rules with a given polynomial accuracy allow the same truncation of a pseudospectral approximation. In practice, however, we care about the cumulative cost of the adaptive algorithm, which must step through successive levels of refinement.

Integration over a bounded interval with uniform weighting offers the widest variety of quadrature choices, and thus allows a thorough comparison. Table \ref{tab:quadCost} summarizes the costs of several common quadrature schemes. First, we see that linear-growth Gaussian quadrature is asymptotically much less efficient than exponential-growth in reaching any particular degree of exactness. However, for rules with fewer than about ten points, this difference is not yet significant.  Second, Clenshaw-Curtis shows efficiency equivalent to exponential-growth Gaussian: both use $n$ points to reach $n$th order polynomial exactness \cite{Clenshaw1960}. However, their performance with respect to external aliasing differs: Clenshaw-Curtis slowly loses accuracy if the integrand is of order greater than $n$, while Gaussian quadrature gives \OO{1} error even on $(n+1)$-order functions \cite{Trefethen2008}. This may make Clenshaw-Curtis Smolyak pseudospectral estimates more efficient. Finally, we consider Gauss-Patterson quadrature, which is nested and has significantly higher polynomial exactness---for a given cumulative cost---than the other types \cite{Patterson1968}. Computing the quadrature points and weights in finite precision (even extended-precision) arithmetic has practically limited Gauss-Patterson rules to 255 points, but we recommend them whenever this is sufficient.
\begin{table}
	\centering
	\footnotesize
		\begin{tabular}{r || c | c | c || c | c | c  || c |c || c |c}
		& \multicolumn{3}{c||}{Lin.\ G} & \multicolumn{3}{c||}{Exp.\ G} & \multicolumn{2}{c||}{C-C} & \multicolumn{2}{c}{G-P} \\
		Order & $p$ & $a$ & $t$ & $p$ & $a$ & $t$ & $p$ & $a$  & $p$ & $a$ \\ \hline
		
		1 & 1 & 1 & 1 & 1 & 1 & 1 & 1 & 1 & 1 & 1  \\
		2 & 2 & 3 & 3 & 2 & 3 & 3 & 3 & 3 & 3 & 5  \\
		3 & 3 & 5 & 6 & 4 & 7 & 7 & 5 & 5 & 7 & 10\\
		4 & 4 & 7 & 10 & 8 & 15 & 15 & 9  & 9 & 15 & 22\\
		5 & 5 & 9 & 15 & 16 & 31 & 31 & 17 & 17 & 31 & 46\\
		6 & 6 & 11 & 21 & 32 & 63 & 63 & 31  & 31 & 63 & 94\\
		$m$ & $m$ & $2m-1$ & $m^2-m/2$ & $2^{m-1}$ & $2^{m}-1 $ & $2^{m}-1 $ & $2^{m-1}+1$ & $2^{m-1}+1$ & &			
		\end{tabular}
	\caption{The cost of four quadrature strategies as their order increases: linear growth Gauss-Legendre quadrature (Lin.\ G), exponential growth Gauss-Legendre quadrature (Exp.\ G), Clenshaw-Curtis quadrature (C-C), and Gauss-Patterson quadrature (G-P). We list the number of points used to compute the given rule (p), the polynomial exactness (a), and the total number of points used so far (t). For nested rule, (p) = (t), so the total column is omitted.}
	\label{tab:quadCost}
\end{table}

For most other weights and intervals, there are fewer choices that provide polynomial exactness, so exponential-growth Gaussian quadrature is our default choice. In the specific case of Gaussian weight, Genz has provided a family of Kronrod extensions, similar to Gauss-Patterson quadrature, which may be a useful option \cite{Genz1996}.

If a linear growth rule is chosen and the domain is symmetric, we suggest that each new level include at least two points, so that the corresponding basis grows by at least one even and one odd basis function. This removes the possibility for unexpected effects on the adaptive strategy if the target function is actually even or odd.

\subsection{Basic convergence: Genz functions}

The Genz family \cite{Genz1984,Genz1987} comprises six parameterized functions, defined from
$[-1,1]^d \to \mathbb{R}$. They are commonly
used to investigate the accuracy of quadrature rules and interpolation
schemes \cite{Barthelmann2000,Klimke2005}. The purpose of this example
is to show that different Smolyak pseudospectral strategies behave
roughly as expected, as evidenced by decreasing $L^2$ approximation
errors as more function evaluations are employed. These functions are as follows: 
\allowdisplaybreaks
\begin{eqnarray*}
\mbox{oscillatory: } f_1(x) &=&  \cos\left(2\pi w_1+\sum_{i=1}^d c_ix_i\right)\\
\mbox{product peak: } f_2(x) &=&  \prod_{i=1}^d\left(c_i^{-2}+(x_i-w_i)^2\right)^{-1}\\
\mbox{corner peak: } f_3(x) &=& \left(1+\sum_{i=1}^d c_ix_i\right)^{-(d+1)}\\
\mbox{Gaussian: } f_4(x) &=& \exp\left(-\sum_{i=1}^d c_i^2\dot(x_i-w_i)^2\right)\\
\mbox{continuous: } f_5(x) &=& \exp\left(-\sum_{i=1}^d c_i^2\dot(|x_i-w_i|)^2\right)\\
\mbox{discontinuous: } f_6(x) &=&  \begin{cases}
	0 & \mbox{if } x_1>w_1 \mbox{ or } x_2> w_2\\
	\exp{\left( \sum_{i=1}^d c_i x_i \right)} & \mbox{otherwise}\\
\end{cases}
\end{eqnarray*} 

Our first test uses five isotropic and \textit{non-adaptive} pseudospectral
approximation strategies. The initial strategy is the isotropic full
tensor pseudospectral algorithm, based on Gauss-Legendre quadrature,
with order growing exponentially with level. The other four strategies
are total-order expansions of increasing order based on the following
quadrature rules: linear growth Gauss-Legendre, exponential growth
Gauss-Legendre, Clenshaw-Curtis, and Gauss-Patterson. All the rules
were selected so that the final rule would have around $10^4$ points.

% For simplicity, all the strategies are non-adaptive and isotropic;
% note that the Genz functions are relatively isotropic

\begin{figure}[htb]
\centering
	\includegraphics{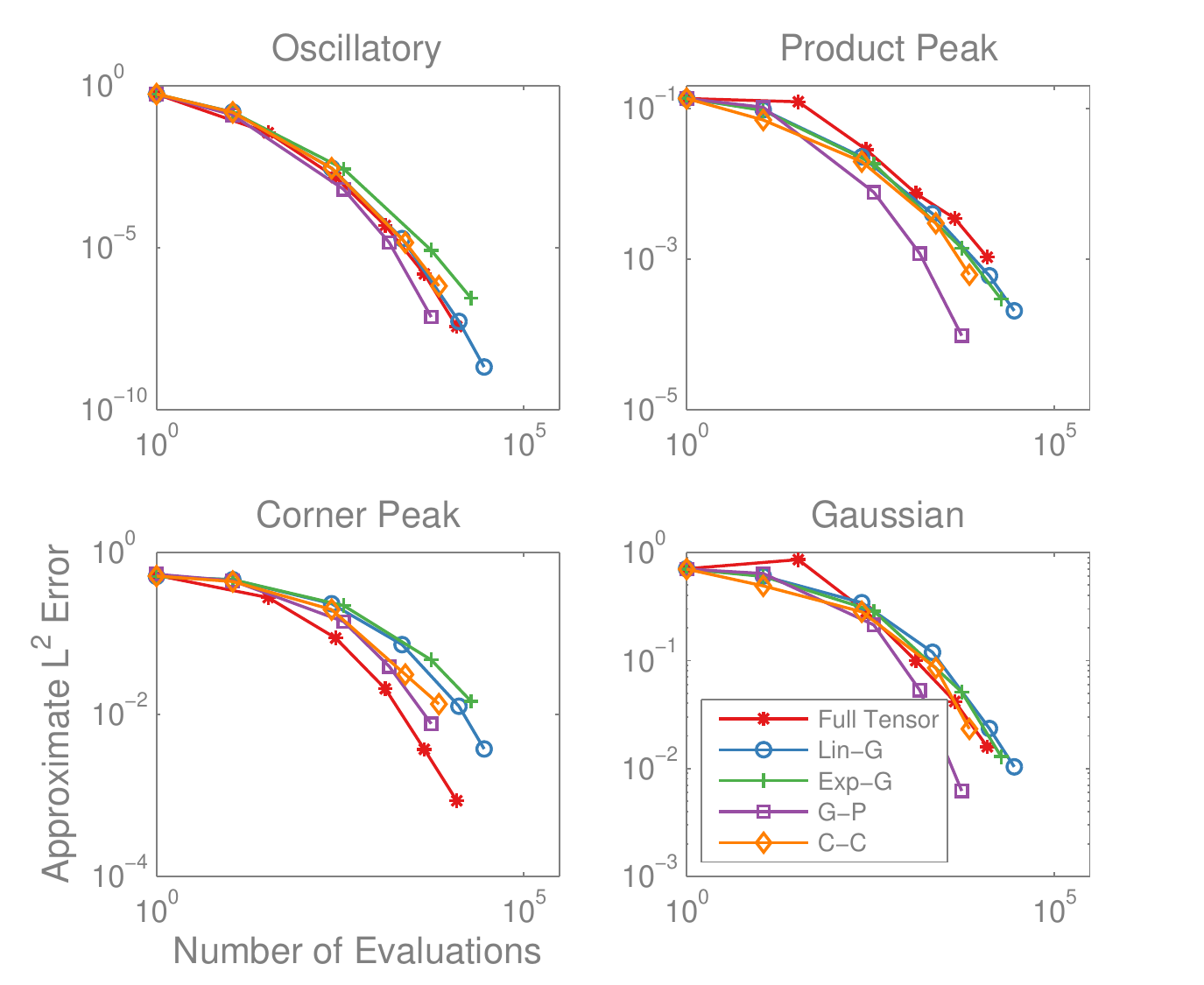}
	\caption{Mean $L^2$ convergence of the non-adaptive isotropic total-order Smolyak pseudospectral algorithm with various quadrature rules, compared to the full tensor pseudospectral algorithm, on the Genz test functions.}
	\label{fig:GenzResults}
\end{figure}

We consider 30 random realizations of each Genz function in $d=5$
dimensions; random parameters for the Genz functions are drawn
uniformly from $[0,1]$, then normalized so that $\|\mathbf{w}\|_1
= 1 $ and $\|\mathbf{c}\|_1 = b_j$, where $j$ indexes the Genz
function type and the constants $b_j$ are as chosen
in~\cite{Barthelmann2000,Klimke2005}.
% \begin{center}
% \begin{tabular}{c|c c c c c c}
% $j$ & 1 &2 &3 &4& 5& 6\\
% \hline
% $b_j$ & 1.5 & $d$ & 1.85& 7.03 & 20.4 & 4.3
% \end{tabular}
% \end{center}
This experiment only uses the first four Genz functions, which are in
$C^\infty$, as pseudospectral methods have well known difficulties on
functions with discontinuities or discontinuous derivatives
\cite{Canuto2006}. Each estimate of $L^2$ approximation error is
computed by Monte Carlo sampling with 10$^4$ samples. Figure
\ref{fig:GenzResults} plots $L^2$ error at each stage, where each
point represents the mean error over the 30 random functions.

Relatively simple conclusions can be drawn from this data. All the
methods show fast convergence, indicating that the internal aliasing issues have
indeed been resolved. In contrast, one would expect direct quadrature
to suffer from large aliasing errors for the three super-linear growth
rules. Otherwise, judging the efficiency of the different rules is not
prudent, because differences in truncation and the structure of the test functions
themselves obscure differences in efficiency. In deference to our
adaptive strategy, we ultimately do not recommend this style of
isotropic and function-independent truncation anyway. 

\begin{figure}[htb]
\centering
	\includegraphics{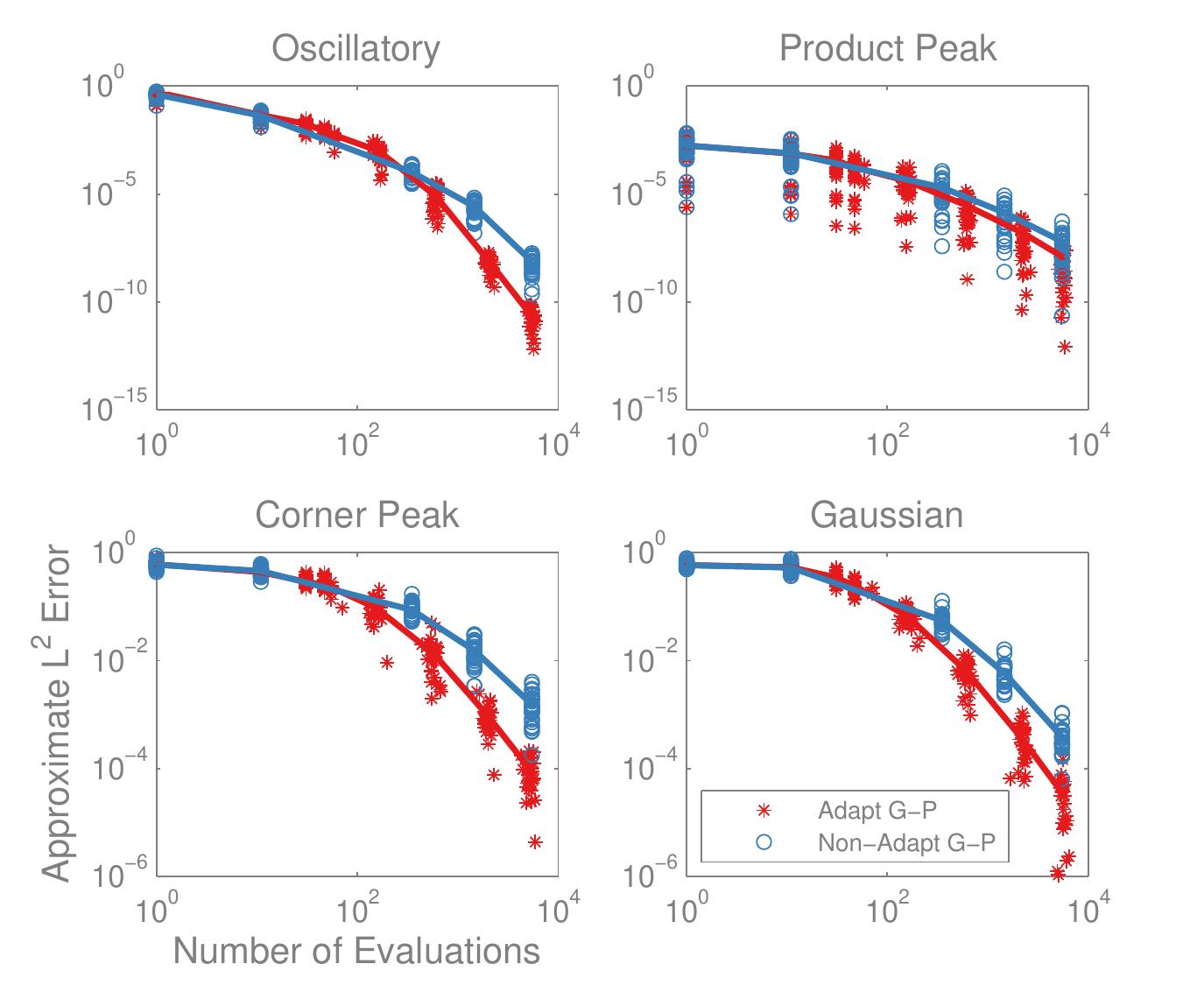}
	\caption{$L^2$ convergence of the adaptive and non-adaptive Gauss-Patterson Smolyak pseudospectral algorithm. Individual results for 30 random instances of the Genz functions are shown.}
	\label{fig:GenzScatter}
\end{figure}

To test our \textit{adaptive} approach, Figure \ref{fig:GenzScatter}
shows results from a similar experiment, now comparing the convergence
of an adaptive Smolyak pseudospectral algorithm with that of a
non-adaptive algorithm. To make the functions less isotropic, we
introduce an exponential decay, replacing each $c_i$ with $c_i
e^{i/5}$, where the $c_i$ are generated and normalized as above. For
consistency, both algorithms are based on Gauss-Patterson
quadrature. As we cannot synchronize the number of evaluations used by
the adaptive algorithm for different functions, we plot individual
errors for the 30 random functions instead of the mean error. This
reveals the variability in difficulty of the functions, which was
hidden in the previous plot. We conclude that the adaptive algorithm
also converges as expected, with performance comparable to or better
than the non-adaptive algorithm. Even though we have included some
anisotropy, these functions include relatively high degrees of
coupling; hence, in this case the non-adaptive strategy is a fairly
suitable choice. For example, the ``product peak'' function shows
little benefit from the adaptive strategy. Although omitted here for
brevity, other quadrature rules produce similar results when comparing
adaptive and non-adaptive algorithms.

\subsection{Adaptivity: chemical kinetics}

To further illustrate the benefits of an adaptive Smolyak approach, we build a
surrogate for a realistic simulation of a combustion kinetics problem.
Specifically, we consider the auto-ignition of a methane-air mixture
given 14 uncertain rate parameters. Governing equations for this
process are a set of stiff nonlinear ODEs expressing conservation of
energy and of chemical species \cite{Kee2003}. The uncertain rate
parameters represent activation energies of reactions governing the conversion of methane to methyl, each endowed with a uniform distribution varying over $[0.8, 1.25]$ of the nominal value. These parameters appear in Arrhenius
expressions for the species production rates, with the reaction
pathways and their nominal rate parameters given by the GRIMech 3.0
mechanism \cite{grimech3:local}. The output
of interest is the logarithm of the ignition time, which is a
functional of the trajectory of the ODE system, and is continuous over the selected parameter ranges. Simulations were
performed with the help of the TChem software library \cite{tchem:local},
which provides convenient evaluations of thermodynamic properties and
species production rates, along with Jacobians for implicit time
integration.

Chemical kinetics are an excellent testbed for adaptive
approximation because, by the nature of detailed kinetic systems, we
expect strong coupling between some inputs and weak coupling between
others, but we cannot predict these couplings \emph{a
  priori}.
We test the effectiveness of adaptive Smolyak pseudospectral methods
based on the four quadrature rules discussed earlier. As our earlier
analysis suggested that Gauss-Patterson quadrature should be most
efficient, our reference solution is a non-adaptive Gauss-Patterson
total-order Smolyak pseudospectral expansion. We ran the non-adaptive
algorithm with a total order index set truncated at $n=5$ (which includes monomial basis terms up through
$\psi_{23}^{(i)}$), using around 40000 point evaluations and taking over
an hour to run. We tuned the four adaptive algorithms to terminate
with approximately the same number of evaluations.

\begin{figure}[htb] \centering
  \includegraphics[scale=.6]{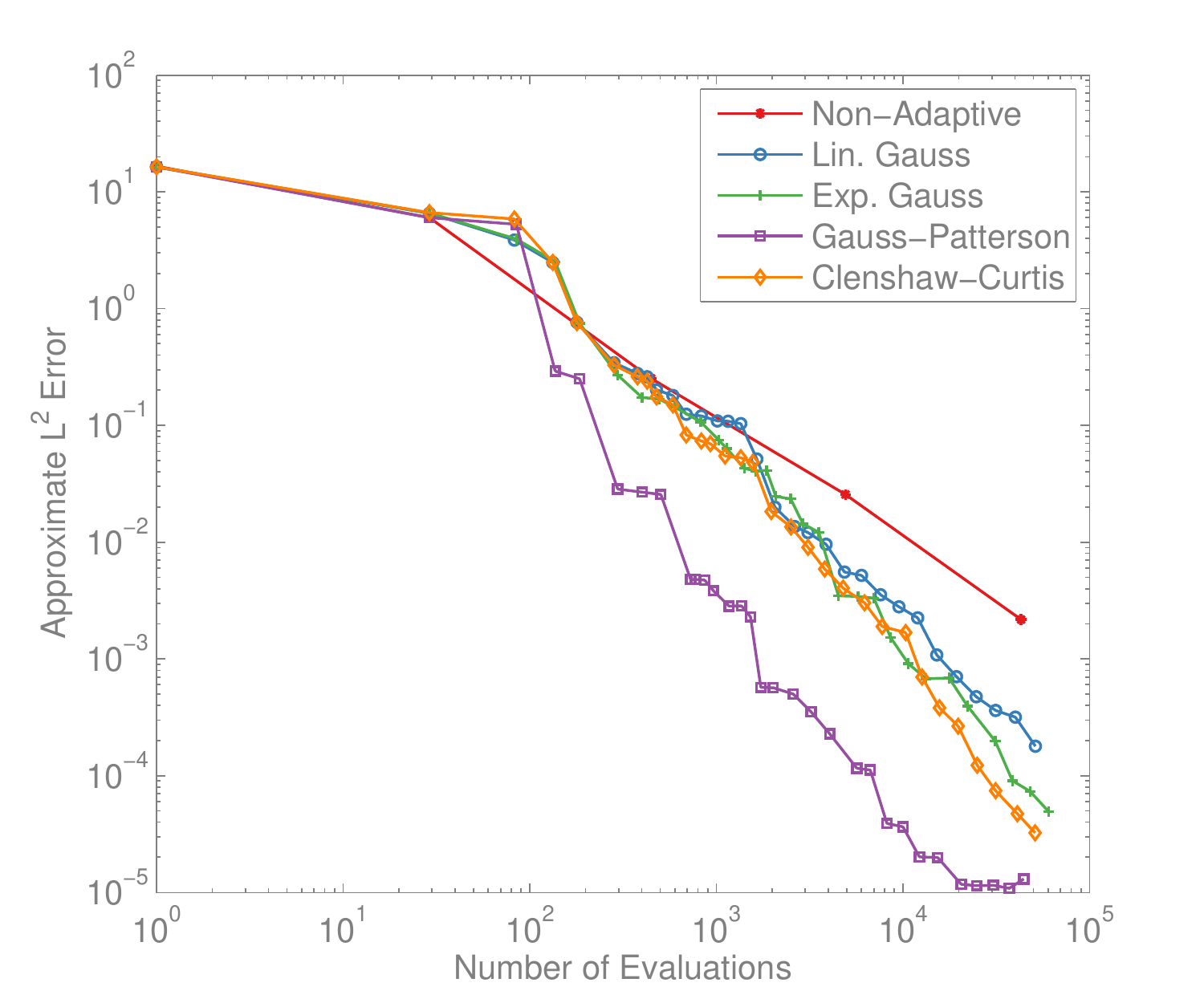}
  \caption{$L^2$ convergence of ignition delay in a 14-dimensional chemical kinetic system; comparing a
    non-adaptive isotropic total-order Gauss-Patterson-based Smolyak
    pseudospectral algorithm to the adaptive algorithm with various
    quadrature rules.}
  \label{fig:combustionConvergence}
\end{figure}

Figure \ref{fig:combustionConvergence} compares convergence of the
five algorithms. The $L^2$ errors reported on the vertical axis are
Monte Carlo estimates using $10^4$ points. Except for a small deviation
at fewer than 200 model evaluations, all of the adaptive methods
significantly outperform the non-adaptive method. The performance of
the different quadrature rules is essentially as predicted in Section
\ref{sec:quadRules}: Gauss-Patterson is the most efficient,
exponential growth Gauss-Legendre and Clenshaw-Curtis are nearly
equivalent, and linear growth Gauss-Legendre performs worse as the
order of the polynomial approximation increases. Compared to the
non-adaptive algorithm, adaptive Gauss-Patterson yields more than two
orders of magnitude reduction in the error at the same number of model
evaluations. Linear growth Gaussian quadrature is initially comparable
to exponential growth Gaussian quadrature, because the asymptotic
benefits of exponential growth do not appear while the algorithm is
principally using very small one-dimensional quadrature rules. At the
end of these experiments, a reasonable number of higher order
quadrature rules are used and the difference becomes visible.

\begin{figure}[htb]	
\centering
	\includegraphics{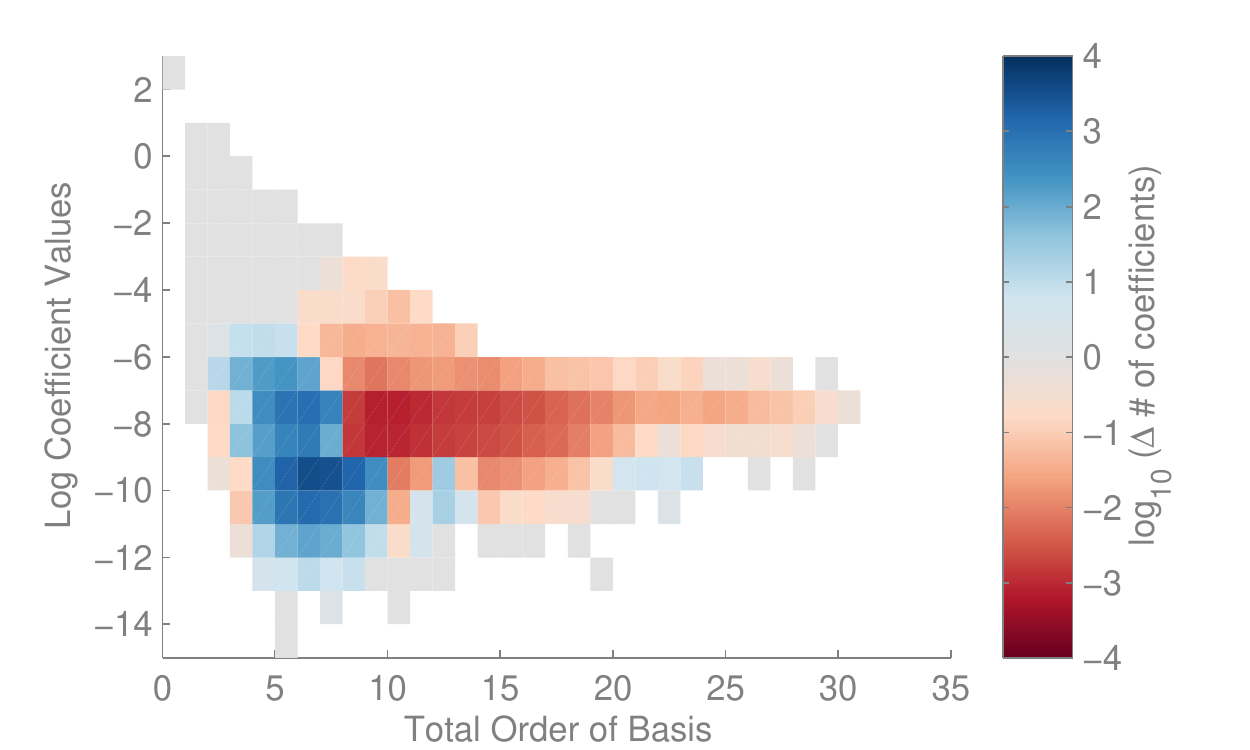}
	\caption{The plot depicts the difference between the
          \emph{number} of coefficients of a particular magnitude and
          order in the final \textit{adaptive} and
          \textit{non-adaptive} Gauss-Patterson based expansions. The
          horizontal axis is the order of the term and the vertical
          axis specifies $\log_{10}$ of the coefficient value. The
          color represents $\log_{10}$ of the difference between
          the two methods, where positive values indicate more terms
          in the non-adaptive expansion. Hence, the dark blue at
          $(6,-10)$ indicates that the non-adaptive expansion includes
          around 3,000 extra terms of magnitude $10^{-10}$ and the
          dark red at $(10,-8)$ indicates that the adaptive expansion
          includes about 1,000 extra terms of magnitude $10^{-8}$.
          Grey squares are the same for both expansions and white
          squares are not present in either.}
	\label{fig:combustionCoeffs}
\end{figure}

We conclude by illustrating that the adaptive algorithm is effective
because it successfully focuses its efforts on high-magnitude
coefficients---that is, coefficients that make the most significant
contributions to the function. Even though the non-adaptive expansion
has around 37,000 terms and the final adaptive Gauss-Patterson
expansion only has about 32,000 terms, the adaptive expansion exhibits
much lower error because most of the additional terms in the
non-adaptive expansion are nearly zero. By skipping many near-zero
coefficients, the adaptive approach is able to locate and estimate a
number of higher-order terms with large magnitudes. Figure
\ref{fig:combustionCoeffs} depicts this pattern by plotting the
difference between the numbers of included terms in the final adaptive
Gauss-Patterson and non-adaptive expansions. The adaptive algorithm
does not actually add any higher order monomials; neither uses
one-dimensional basis terms of order higher than $\psi^{(i)}_{23}$.
Instead, the adaptive algorithm adds mixed terms of higher total
order, thus capturing the coupling of certain variables in more detail
than the non-adaptive algorithm. The figure shows that terms through
30\textsuperscript{th} order are included in the adaptive expansion,
all of which are products of non-constant polynomials in more than one
dimension.

\subsection{Performance of the global error indicator}

\begin{figure}[htb] \centering
  \includegraphics[scale=.65]{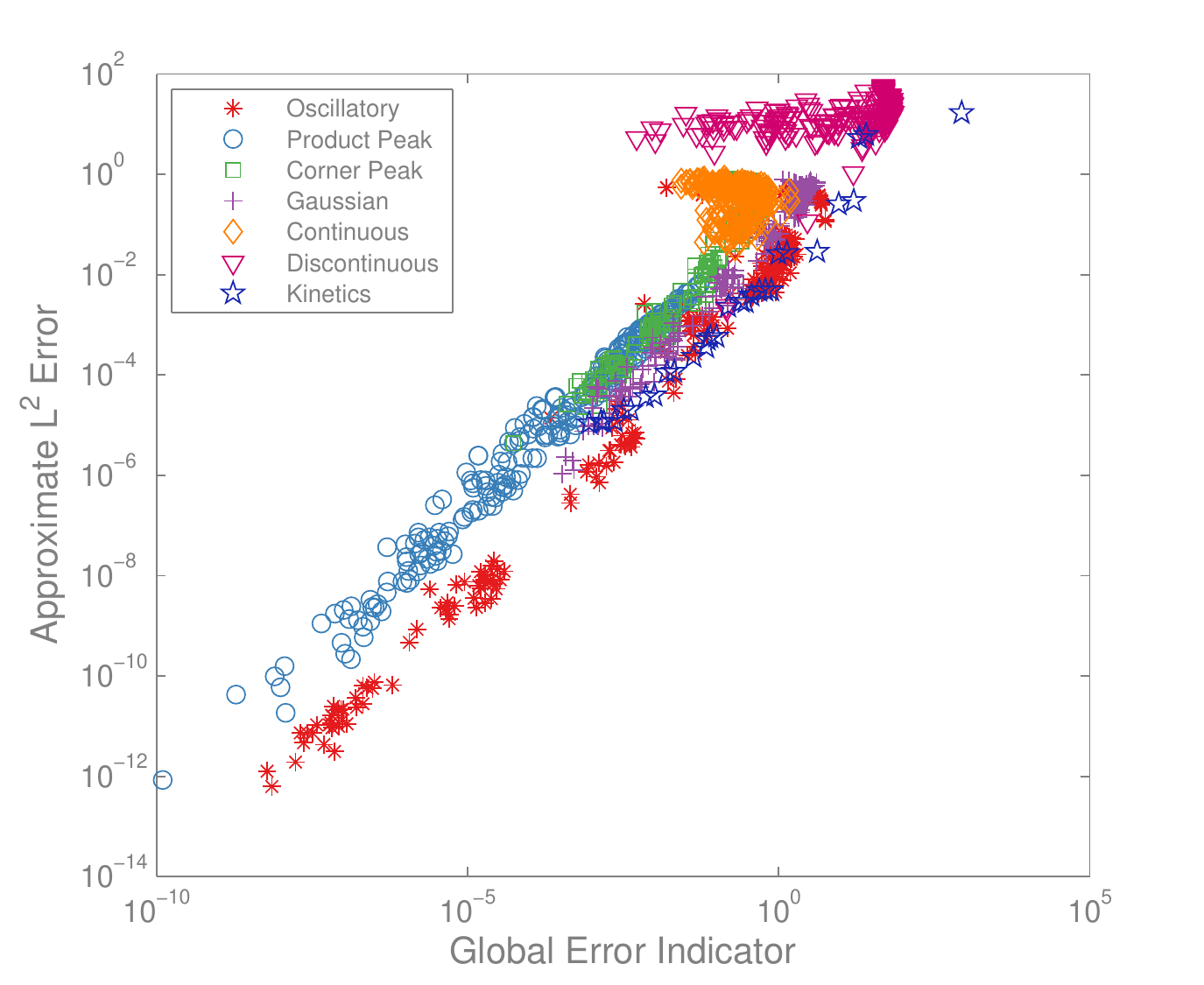}
  \caption{Relationship between the termination criterion
    (\ref{eq:globalError}) and the estimated $L^2$ error for every
    function tested.}
  \label{fig:terminationPlot}
\end{figure}

To evaluate the termination criterion, we collected the global error
indicator during runs of the adaptive algorithm for all of the test functions
described above, including the slowly converging non-smooth Genz
functions omitted before. The discontinuous Genz function does not
include the exponential coefficient decay because the discontinuity
already creates strong anisotropy. Results are shown for
Gauss-Patterson quadrature. The relationship between the estimated
$L^2$ error and the global error indicator $\epsilon_g$ is shown in
Figure \ref{fig:terminationPlot}. For the smooth test functions,
$\epsilon_g$ is actually an excellent indicator, as it is largely
within an order of magnitude of the correct value and essentially linearly related to it. However, the non-smooth Genz functions illustrate the hazard of relying too heavily on this indicator: although the adaptive algorithm does decrease both the errors and the indicator, the relationship between the two appears far less direct.

\section{Conclusions}

This paper gives a rigorous development of Smolyak pseudospectral algorithms, a practical approach for constructing polynomial chaos expansions from point evaluations of a function. A common alternative approach, direct quadrature, has previously been shown to suffer from large errors. We explain these errors as a consequence of internal aliasing and delineate the exact circumstances, derived from properties of the chosen polynomial basis and quadrature rules, under which internal aliasing will occur. Internal aliasing is a problem inherent to direct quadrature approaches, which use a single (sparse) quadrature rule to compute a set of spectral coefficients.
% Here: In a broader sense, what was wrong with this algorithm? Never mind synchronization issues, which are discussed in the next few lines.
%
These approaches fail because they substitute a numerical approximation for only a portion of the algorithm, i.e., the evaluation of integrals, without considering the impact of this approximation on the entire construction. For almost all sparse quadrature rules, internal aliasing errors may be overcome only through an inefficient use of function evaluations.
In contrast, the Smolyak pseudospectral algorithm computes spectral coefficients by assembling tensor-product pseudospectral approximations in a coherent fashion that avoids internal aliasing by construction; moreover, it has smaller external aliasing errors. To establish these properties, we extend the known result that the exact set of a Smolyak pseudospectral approximation contains a union of the exact sets of all its constituent tensor-product approximation operators to the case of arbitrary admissible Smolyak multi-index sets. These results are applicable to any choice of quadrature rule and generalized sparse grid, and are verified through numerical demonstrations; hence, we suggest that the Smolyak pseudospectral algorithm is a superior approach in almost all contexts.

A key strength of Smolyak algorithms is that they are highly customizable through the choice of admissible multi-index sets. To this end, we describe a simple alteration to the adaptive sparse quadrature approaches of \cite{Gerstner2003, Hegland2003}, creating a corresponding method for adaptive pseudospectral approximation. Numerical experiments then evaluate the performance of different quadrature rules and of adaptive versus non-adaptive pseudospectral approximation. Tests of the adaptive method on a realistic chemical kinetics problem show multiple order-of-magnitude gains in accuracy over a non-adaptive approach. Although the adaptive strategy will not improve approximation performance for every function, we have little evidence that it is ever harmful and hence widely recommend its use.

While the adaptive approach illustrated here is deliberately simple, many extensions are possible. For instance, as described in Section~\ref{s:adaptivecomments}, measures of computational cost may be added to the dimension refinement criterion. One could also use the gradient of the $L^2$ error indicator to identify optimal directions in the space of multi-indices along which to continue refinement, or to avoid adding all the forward neighbors of the multi-index selected for refinement. These and other developments will be pursued in future work.

A flexible open-source \CC \ code implementing the adaptive approximation method discussed in this paper is available at \url{https://bitbucket.org/mituq/muq/}.

%%%%%%%%%%%%%%%%%%%%%%%%%%%%%%%%%%%%%%%%%%%%%%%%%

%%%%%%%%%%%%%%%%%%%

% We have given a theoretical foundation for the observed performance problems with direct quadrature, explaining that these errors are an unavoidable consequence of aliasing, as derived from properties of the polynomial basis and quadrature rules. In contrast, the Smolyak pseudospectral algorithm does not suffer from internal aliasing, and we provided well known examples to illustrate their convergence. Furthermore, Smolyak algorithms, if we use generalized Smolyak multi-index sets, are greatly customizable. Finally, we described a simple alteration to the successful adaptive algorithm for sparse quadrature of Gerstner and Griebel, and showed that the adaptive method performs well in practice on a realistic chemical kinetics problem \cite{Gerstner2003}. 

% We believe that these new results for efficiently building PCEs for loosely coupled problems should allow surrogate based uncertainty quantification methods to become more widely applicable and useful.

% Delineating circumstances under which aliasing errors will or will not occur. While some (e.g., linear-growth GQ) might not, 

% quite general: not limited to Gaussian quadrature or nested quad rules, any admissible sparse grid; 

% evaluated performance of different quadrature rules, with particularly favorable performance of gauss-patterson

% extension =  additional modes of refinement: cost, gradient, directionality, etc

\section*{Acknowledgments}
The authors would like to thank Paul Constantine for helpful
discussions and for sharing a preprint of his paper, which inspired
this work. We would also like to thank Omar Knio and Justin Winokur
for many helpful discussions, and Tom Coles for help with the chemical
kinetics example. P.\ Conrad was supported during this work by a
Department of Defense NDSEG Fellowship and an NSF graduate
fellowship. P.\ Conrad and Y.\ Marzouk acknowledge additional support
from the Scientific Discovery through Advanced Computing (SciDAC)
program funded by the US Department of Energy, Office of Science,
Advanced Scientific Computing Research under award number
DE-SC0007099.

%%%%%%%%%%%%%%%%%%%%%%%%%%%%%%%%%%%%%%%%%%%%%%%%%%%%%%%%%%%%%
%% ==> Some hints are following:

%%%%%%%%%%%%%%%%%%%%%%%%%%%%%%%%%%%%%%%%%%%%%%%%%%%%%%%%%%%%%
%% BIBLIOGRAPHY AND OTHER LISTS
%%%%%%%%%%%%%%%%%%%%%%%%%%%%%%%%%%%%%%%%%%%%%%%%%%%%%%%%%%%%%
%% A small distance to the other stuff in the table of contents (toc)
%\addtocontents{toc}{\protect\vspace*{\baselineskip}}

%% The Bibliography
%% ==> You need a file 'literature.bib' for this.
%% ==> You need to run BibTeX for this (Project | Properties... | Uses BibTeX)
%\addcontentsline{toc}{chapter}{Bibliography} %'Bibliography' into toc
%\nocite{*} %Even non-cited BibTeX-Entries will be shown.
\bibliographystyle{siam} %Style of Bibliography: plain / apalike / amsalpha / ...
\bibliography{local,library} %You need a file 'literature.bib' for this.

%% The List of Figures
%\clearpage
%\addcontentsline{toc}{chapter}{List of Figures}
%\listoffigures

%% The List of Tables
%\clearpage
%\addcontentsline{toc}{chapter}{List of Tables}
%\listoftables

%%%%%%%%%%%%%%%%%%%%%%%%%%%%%%%%%%%%%%%%%%%%%%%%%%%%%%%%%%%%%
%% APPENDICES
%%%%%%%%%%%%%%%%%%%%%%%%%%%%%%%%%%%%%%%%%%%%%%%%%%%%%%%%%%%%%
\appendix

%% ==> Write your text here or include other files.

%\input{FileName} %You need a file 'FileName.tex' for this.

\end{document}